\let\myfrac=\frac
\input eplain
\let\frac=\myfrac
\input amstex
\input epsf




\loadeufm \loadmsam \loadmsbm
\message{symbol names}\UseAMSsymbols\message{,}

\font\myfontdefault=cmr10

\font\mytdmchapfont=cmb10 at 14pt
\font\mytdmheadfont=cmb10 at 10pt
\font\mytdmsubheadfont=cmr10

\magnification 1200
\newif\ifinappendices
\newif\ifundefinedreferences
\newif\ifchangedreferences
\newif\ifloadreferences
\newif\ifmakebiblio
\newif\ifmaketdm

\undefinedreferencestrue
\changedreferencesfalse


\loadreferencestrue
\makebibliofalse
\maketdmfalse

\def\headpenalty{-400}     
\def\proclaimpenalty{-200} 

%
%

\def\alphanum#1{\ifcase #1 _\or A\or B\or C\or D\or E\or F\or G\or H\or I\or J\or K\or L\or M\or N\or O\or P\or Q\or R\or S\or T\or U\or V\or W\or X\or Y\or Z\fi}
\def\gobbleeight#1#2#3#4#5#6#7#8{}

\newwrite\references
\newwrite\tdm
\newwrite\biblio

\newcount\chapno
\newcount\headno
\newcount\subheadno
\newcount\procno
\newcount\figno
\newcount\citationno

\def\setcatcodes{%
\catcode`\!=0 \catcode`\\=11}%

\ifloadreferences
    {\catcode`\@=11 \catcode`\_=11%
    \global\def\_@citation@CaffI{1}
\global\def\_@citation@CaffII{2}
\global\def\_@citation@CaffNirSprI{3}
\global\def\_@citation@CaffNirSprII{4}
\global\def\_@citation@CaffNirSprV{5}
\global\def\_@citation@CalI{6}
\global\def\_@citation@Gauss{7}
\global\def\_@citation@GilbTrud{8}
\global\def\_@citation@Guan{9}
\global\def\_@citation@GuanSpruckO{10}
\global\def\_@citation@GuanSpruckI{11}
\global\def\_@citation@GuanSpruckII{12}
\global\def\_@citation@GuanSpruckIII{13}
\global\def\_@citation@Gut{14}
\global\def\_@citation@Jorgens{15}
\global\def\_@citation@LabI{16}
\global\def\_@citation@LoftI{17}
\global\def\_@citation@LoftII{18}
\global\def\_@citation@Pog{19}
\global\def\_@citation@PogB{20}
\global\def\_@citation@RosSpruck{21}
\global\def\_@citation@SmiPPG{22}
\global\def\_@citation@SmiPPH{23}
\global\def\_@citation@SmiNLD{24}
\global\def\_@citation@SmiFCS{25}
\global\def\_@citation@SmiSLC{26}
\global\def\_@citation@Spruck{27}
\global\def\_@citation@TrudWang{28}
\global\def\_@proc@TheoremExistence{1.1}
\global\def\_@proc@TheoremBVP{1.2}
\global\def\_@proc@LemmaFirstOrderBounds{2.1}
\global\def\_@head@HeadGaussCurvature{3}
\global\def\_@proc@PropFormulaForGaussCurvature{3.1}
\global\def\_@proc@PropConnexionInGraph{3.2}
\global\def\_@head@HeadMaximumPrinciples{4}
\global\def\_@head@HeadSecondOrderBoundaryBounds{5}
\global\def\_@proc@PropSecondOrderBoundaryBounds{5.1}
\global\def\_@proc@PropFirstBarrierEstimate{5.3}
\global\def\_@proc@CorSecondBarrierEstimate{5.5}
\global\def\_@proc@PropThirdBarrierEstimate{5.6}
\global\def\_@head@HeadSecondOrderBoundsOverTheInterior{6}
\global\def\_@proc@PropSecondOrderBounds{6.1}
\global\def\_@proc@LemmaCommutationRelations{6.3}
\global\def\_@proc@CorSecondDerivativesOfA{6.4}
\global\def\_@proc@PropSuperharmonic{6.6}
\global\def\_@proc@CorSuperharmonic{6.7}
\global\def\_@proc@CorMaximumPrinciple{6.9}
\global\def\_@head@HeadCompactness{7}
\global\def\_@proc@LemmaKisProper{7.1}
\global\def\_@proc@CorKisProper{7.2}
\global\def\_@head@HeadOneDimensionalFamiliesOfSolutions{8}
\global\def\_@proc@PropTransversality{8.1}
\global\def\_@proc@PropCompactnessOfGamma{8.2}
\global\def\_@head@HeadLocalAndGlobalRigidity{9}
\global\def\_@proc@DefnLocalAndGlobalRigidity{9.1}
\global\def\_@proc@PropPropertiesOfRigidity{9.2}
\global\def\_@proc@LemmaExistence{10.2}
\global\def\_@head@HeadHyperbolicCase{11}
\global\def\_@proc@LemmaClosedMeansInConvexHull{11.1}
\global\def\_@proc@PropNodalSetIsClosed{11.2}
\global\def\_@proc@LemmaDirichletInHyperbolicSpace{11.3}
\global\def\_@head@HeadRelationsToExistingResults{12}
\global\def\_@proc@TheoremGuan{12.1}
\global\def\_@proc@TheoremRosSpruck{12.2}
\global\def\_@proc@LemmaSmoothnessOfLimit{A.1}
\global\def\_@proc@LemmaBasicSquareRelation{A.2}
    }%
\else
    \openout\references=references.tex
\fi

\newcount\newchapflag 
\newcount\showpagenumflag 

\global\chapno = -1 
\global\citationno=0
\global\headno = 0
\global\subheadno = 0
\global\procno = 0
\global\figno = 0

\def\resetcounters{%
\global\headno = 0%
\global\subheadno = 0%
\global\procno = 0%
\global\figno = 0%
}

\global\newchapflag=0 
\global\showpagenumflag=0 

\def\chinfo{\ifinappendices\alphanum\chapno\else\the\chapno\fi}%
\def\headinfo{\ifinappendices\alphanum\headno\else\the\headno\fi}%
\def\subheadinfo{\headinfo.\the\subheadno}
\def\procinfo{\headinfo.\the\procno}
\def\figinfo{\the\figno}        
\def\citationinfo{\the\citationno}%
\def\nextheadno{\global\advance\headno by 1 \global\subheadno = 0 \global\procno = 0}
\def\nextsubheadno{\global\advance\subheadno by 1}
\def\nextprocno{\global\advance\procno by 1 \procinfo}
\def\nextfigno{\global\advance\figno by 1 \figinfo}

{\global\let\noe=\noexpand%
%
%
\catcode`\@=11%
\catcode`\_=11%
\setcatcodes%
!global!def!_@@internal@@makeref#1{%
!global!expandafter!def!csname #1ref!endcsname##1{%
!csname _@#1@##1!endcsname%
!expandafter!ifx!csname _@#1@##1!endcsname!relax%
    !write16{#1 ##1 not defined - run saving references}%
    !undefinedreferencestrue%
!fi}}%
!global!def!_@@internal@@makelabel#1{%
!global!expandafter!def!csname #1label!endcsname##1{%
!edef!temptoken{!csname #1info!endcsname}%
!ifloadreferences%
    !expandafter!ifx!csname _@#1@##1!endcsname!relax%
        !write16{#1 ##1 not hitherto defined - rerun saving references}%
        !changedreferencestrue%
    !else%
        !expandafter!ifx!csname _@#1@##1!endcsname!temptoken%
        !else
            !write16{#1 ##1 reference has changed - rerun saving references}%
            !changedreferencestrue%
        !fi%
    !fi%
!else%
    !expandafter!edef!csname _@#1@##1!endcsname{!temptoken}%
    !edef!textoutput{!write!references{\global\def\_@#1@##1{!temptoken}}}%
    !textoutput%
!fi}}%
!global!def!makecounter#1{!_@@internal@@makelabel{#1}!_@@internal@@makeref{#1}}%
!unsetcatcodes%
}
\makecounter{ch}%
\makecounter{head}%
\makecounter{subhead}%
\makecounter{proc}%
\makecounter{fig}%
\makecounter{citation}%
\def\newref#1#2{%
\def\temptext{#2}%
\edef\bibliotextoutput{\expandafter\gobbleeight\meaning\temptext}%
\global\advance\citationno by 1\citationlabel{#1}%
\ifmakebiblio%
    \edef\fileoutput{\write\biblio{\noindent\hbox to 0pt{\hss$[\the\citationno]$}\hskip 0.2em\bibliotextoutput\medskip}}%
    \fileoutput%
\fi}%
\def\cite#1{%
$[\citationref{#1}]$%
\ifmakebiblio%
    \edef\fileoutput{\write\biblio{#1}}%
    \fileoutput%
\fi%
}%
%
%
%

\let\mypar=\par


\def\raggedleft{\leftskip=0pt plus 1fil \parfillskip=0pt}


\font\lettrinefont=cmr10 at 28pt
\def\lettrine #1[#2][#3]#4%
{\hangafter -#1 \hangindent #2
\noindent\hskip -#2 \vtop to 0pt{
\kern #3 \hbox to #2 {\lettrinefont #4\hss}\vss}}

\font\mylettrinefont=cmr10 at 28pt
\def\mylettrine #1[#2][#3][#4]#5%
{\hangafter -#1 \hangindent #2
\noindent\hskip -#2 \vtop to 0pt{
\kern #3 \hbox to #2 {\mylettrinefont #5\hss}\vss}}


\edef\Pagetitle={Blank}

\headline={\hfil\Pagetitle\hfil}

\footline={\hfil\myfontdefault\folio\hfil}

\def\nextoddpage
{
\newpage%
\ifodd\pageno%
\else%
    \global\showpagenumflag = 0%
    \null%
    \vfil%
    \eject%
    \global\showpagenumflag = 1%
\fi%
}


\def\newchap#1#2%
{%
%
%
\global\advance\chapno by 1%
\resetcounters%
%
%
\newpage%
\ifodd\pageno%
\else%
    \global\showpagenumflag = 0%
    \null%
    \vfil%
    \eject%
    \global\showpagenumflag = 1%
\fi%
\global\newchapflag = 1%
\global\showpagenumflag = 1%
%
%
{\font\chapfontA=cmsl10 at 30pt%
\font\chapfontB=cmsl10 at 25pt%
\null\vskip 5cm%
{\chapfontA\raggedleft\hfil%
{%
\ifnum\chapno=0
    \phantom{%
    \ifinappendices%
        Annexe \alphanum\chapno%
    \else%
        \the\chapno%
    \fi}%
\else%
    \ifinappendices%
        Annexe \alphanum\chapno%
    \else%
        \the\chapno%
    \fi%
\fi%
}%
\par}%
\vskip 2cm%
{\chapfontB\raggedleft%
\lineskiplimit=0pt%
\lineskip=0.8ex%
\hfil #1\par}%
\vskip 2cm%
}%
\edef\Pagetitle{#2}%
%
%
\ifmaketdm%
    \def\temp{#2}%
    \def\tempbis{\nobreak}%
    \edef\chaptitle{\expandafter\gobbleeight\meaning\temp}%
    \edef\mynobreak{\expandafter\gobbleeight\meaning\tempbis}%
    \edef\textoutput{\write\tdm{\bigskip{\noexpand\mytdmchapfont\noindent\chinfo\ - \chaptitle\hfill\noexpand\folio}\par\mynobreak}}%
\fi%
\textoutput%
}


\def\newhead#1%
{%
\ifhmode%
    \mypar%
\fi%
\ifnum\headno=0%
\ifinappendices
    \nobreak\vskip -\lastskip%
    \nobreak\vskip .5cm%
\fi
\else%
    \nobreak\vskip -\lastskip%
    \nobreak\vskip .5cm%
\fi%
\nextheadno%
\ifmaketdm%
    \def\temp{#1}%
    \edef\sectiontitle{\expandafter\gobbleeight\meaning\temp}%
    \edef\textoutput{\write\tdm{\noindent{\noexpand\mytdmheadfont\quad\headinfo\ - \sectiontitle\hfill\noexpand\folio}\par}}%
    \textoutput%
\fi%
\font\headfontA=cmbx10 at 14pt%
{\headfontA\noindent\headinfo\ - #1.\hfil}%
\nobreak\vskip .5cm%
}%


\def\newsubhead#1%
{%
\ifhmode%
    \mypar%
\fi%
\ifnum\subheadno=0%
\else%
    \penalty\headpenalty\vskip .4cm%
\fi%
\nextsubheadno%
\ifmaketdm%
    \def\temp{#1}%
    \edef\subsectiontitle{\expandafter\gobbleeight\meaning\temp}%
    \edef\textoutput{\write\tdm{\noindent{\noexpand\mytdmsubheadfont\quad\quad\subheadinfo\ - \subsectiontitle\hfill\noexpand\folio}\par}}%
    \textoutput%
\fi%
\font\subheadfontA=cmsl10 at 12pt
{\subheadfontA\noindent\subheadinfo\ #1.\hfil}%
\nobreak\vskip .25cm %
}%

%
%


\font\mathromanten=cmr10
\font\mathromanseven=cmr7
\font\mathromanfive=cmr5
\newfam\mathromanfam
\textfont\mathromanfam=\mathromanten
\scriptfont\mathromanfam=\mathromanseven
\scriptscriptfont\mathromanfam=\mathromanfive
\def\mathroman{\fam\mathromanfam}


\font\sf=cmss12

\font\sansseriften=cmss10
\font\sansserifseven=cmss7
\font\sansseriffive=cmss5
\newfam\sansseriffam
\textfont\sansseriffam=\sansseriften
\scriptfont\sansseriffam=\sansserifseven
\scriptscriptfont\sansseriffam=\sansseriffive
\def\mathsf{\fam\sansseriffam}


\font\bftwelve=cmb12

\font\boldten=cmb10
\font\boldseven=cmb7
\font\boldfive=cmb5
\newfam\mathboldfam
\textfont\mathboldfam=\boldten
\scriptfont\mathboldfam=\boldseven
\scriptscriptfont\mathboldfam=\boldfive
\def\mathbf{\fam\mathboldfam}


\font\mycmmiten=cmmi10
\font\mycmmiseven=cmmi7
\font\mycmmifive=cmmi5
\newfam\mycmmifam
\textfont\mycmmifam=\mycmmiten
\scriptfont\mycmmifam=\mycmmiseven
\scriptscriptfont\mycmmifam=\mycmmifive

\def\hexa#1{\ifcase #1 0\or 1\or 2\or 3\or 4\or 5\or 6\or 7\or 8\or 9\or A\or B\or C\or D\or E\or F\fi}
\mathchardef\mathi="7\hexa\mycmmifam7B
\mathchardef\mathj="7\hexa\mycmmifam7C


\font\mymsbmten=msbm10 at 8pt
\font\mymsbmseven=msbm7 at 5.6pt
\font\mymsbmfive=msbm5 at 4pt
\newfam\mymsbmfam
\textfont\mymsbmfam=\mymsbmten
\scriptfont\mymsbmfam=\mymsbmseven
\scriptscriptfont\mymsbmfam=\mymsbmfive

\mathchardef\mybeth="7\hexa\mymsbmfam69
\mathchardef\mygimmel="7\hexa\mymsbmfam6A
\mathchardef\mydaleth="7\hexa\mymsbmfam6B


\def\placelabel[#1][#2]#3{{%
\setbox10=\hbox{\raise #2cm \hbox{\hskip #1cm #3}}%
\ht10=0pt%
\dp10=0pt%
\wd10=0pt%
\box10}}%


\newif\ifinproclaim%
\global\inproclaimfalse%
\def\proclaim#1{%
\medskip%
%
%
\bgroup%
\inproclaimtrue%
\setbox10=\vbox\bgroup\leftskip=0.8em\noindent{\bftwelve #1}\sf%
}

\def\endproclaim{%
\egroup%
\setbox11=\vtop{\noindent\vrule height \ht10 depth \dp10 width 0.1em}%
\wd11=0pt%
\setbox12=\hbox{\copy11\kern 0.3em\copy11\kern 0.3em}%
\wd12=0pt%
\setbox13=\hbox{\noindent\box12\box10}%
\noindent\unhbox13%
\egroup%
\medskip\ignorespaces%
}

\def\proclaim#1{%
\medskip%
\bgroup%
\inproclaimtrue%
\noindent{\bftwelve #1}%
\nobreak\medskip%
\sf%
}

\def\endproclaim{%
\mypar\egroup\penalty\proclaimpenalty\medskip\ignorespaces%
}

\def\noskipproclaim#1{%
\medskip%
\bgroup%
\inproclaimtrue%
\noindent{\bf #1}\nobreak\sl%
}

\def\endnoskipproclaim{%
\mypar\egroup\penalty\proclaimpenalty\medskip\ignorespaces%
}


\def\ninn{{n\in\Bbb{N}}}
\def\minn{{m\in\Bbb{N}}}

\def\proof{{\noindent\bf Proof:\ }}

\def\remark{{\noindent\sl Remark:\ }}

\def\msf#1{{\mathsf #1}}

\def\qed{~$\square$}
\def\munion{\mathop{\cup}}
\def\minter{\mathop{\cap}}
\def\myitem#1{%
    \noindent\hbox to .5cm{\hfill#1\hss}
}

\catcode`\@=11
\def\Eqalign#1{\null\,\vcenter{\openup\jot\m@th\ialign{%
\strut\hfil$\displaystyle{##}$&$\displaystyle{{}##}$\hfil%
&&\quad\strut\hfil$\displaystyle{##}$&$\displaystyle{{}##}$%
\hfil\crcr #1\crcr}}\,}
\catcode`\@=12

\def\makeop#1{%
\global\expandafter\def\csname op#1\endcsname{{\mathroman #1}}}%

\def\makeopsmall#1{%
\global\expandafter\def\csname op#1\endcsname{{\mathroman{\lowercase{#1}}}}}%

\makeopsmall{ArcTan}%
\makeopsmall{ArcCos}%
\makeop{Arg}%
\makeop{Det}%
\makeop{Log}%
\makeop{Re}%
\makeop{Im}%
\makeop{Dim}%
\makeopsmall{Tan}%
\makeop{Ker}%
\makeopsmall{Cos}%
\makeopsmall{Sin}%
\makeop{Exp}%
\makeopsmall{Tanh}%
\makeop{Tr}%
\makeop{End}%
\makeop{Long}%
\makeop{Ch}%
\makeop{Exp}%
\makeop{Eval}%
\makeop{Lift}%
\makeop{Int}%
\makeop{Ext}%
\makeop{Aire}%
\makeop{Im}%
\makeop{Conf}%
\makeop{Exp}%
\makeop{Mod}%
\makeop{Log}%
\makeop{Sgn}%
\makeop{Ext}%
\makeop{Int}%
\makeop{Ln}%
\makeop{Dist}%
\makeop{Aut}%
\makeop{Id}%
\makeop{GL}%
\makeop{SO}%
\makeop{Homeo}%
\makeop{Vol}%
\makeop{Ric}%
\makeop{Hess}%
\makeop{Euc}%
\makeop{Isom}%
\makeop{Max}%
\makeop{SW}%
\makeop{SL}%
\makeop{Long}%
\makeop{Fix}%
\makeop{Wind}%
\makeop{Diag}%
\makeop{dVol}%
\makeop{Symm}%
\makeop{Ad}%
\makeop{Diam}%
\makeop{loc}%
\makeopsmall{Sinh}%
\makeopsmall{Cosh}%
\makeop{Len}%
\makeop{Length}%
\makeop{Conv}%
\makeop{Min}%
\makeop{Area}%
\font\mycirclefont=cmsy7
\def\textcircle{{\raise 0.3ex \hbox{\mycirclefont\char'015}}}

\let\emph=\bf

\hyphenation{quasi-con-formal}

%
%

\ifmakebiblio%
    \openout\biblio=biblio.tex %
    {%
        \edef\fileoutput{\write\biblio{\bgroup\leftskip=2em}}%
        \fileoutput
    }%
\fi%

\newref{CaffI}{Caffarelli L., A localization property of viscosity solutions to the Monge-Amp\`ere equation and their strict convexity, {\sl Annals of Math.} {\bf 131} (1990), 129---134}
\newref{CaffII}{Caffarelli L., {\sl Monge-Amp\'ere equation, div-curl theorems in Lagrangian coordinates, compression and rotation}, Lecture Notes, 1997}
\newref{CaffNirSprI}{Caffarelli L., Nirenberg L., Spruck J., The Dirichlet problem for nonlinear second-Order elliptic equations. I. Monge Amp\`ere equation, {\sl Comm. Pure Appl. Math.} {\bf 37} (1984), no. 3, 369--402}
\newref{CaffNirSprII}{Caffarelli L., Kohn J. J., Nirenberg L., Spruck J., The Dirichlet problem for nonlinear second-order elliptic equations. II. Complex Monge Amp\`ere, and uniformly elliptic, equations, {\sl Comm. Pure Appl. Math.} {\bf 38} (1985), no. 2, 209--252}
\newref{CaffNirSprV}{Caffarelli L., Nirenberg L., Spruck J., Nonlinear second-order elliptic equations. V. The Dirichlet problem for Weingarten hypersurfaces, {\sl Comm. Pure Appl. Math.} {\bf 41} (1988), no. 1, 47--70}
\newref{CalI}{Calabi E., Improper affine hypersurfaces of convex type and a generalization of a theorem by K. J\"orgens, {\sl Mich. Math. J.} {\bf 5} (1958), 105--126}
\newref{Gauss}{Gauss C. F., General investigations of curved surfaces, Translated from the Latin and German by Adam Hiltebeitel and James Morehead, Raven Press, Hewlett, N.Y.}
\newref{GilbTrud}{Gilbarg D., Trudinger N. S., {\sl Elliptic Partial Differential Equations of Second Order}, Classics in Mathematics, Springer, Berlin, (2001)}
\newref{Guan}{Guan B., The Dirichlet problem for Monge-Amp\`ere equations in non-convex domains and spacelike hypersurfaces of constant Gauss curvature, {\sl Trans. Amer. Math. Soc.} {\bf 350} (1998), 4955--4971}
\newref{GuanSpruckO}{Guan B., Spruck J., Boundary value problems on $\Bbb{S}^n$ for surfaces of constant Gauss curvature, {\sl Ann. of Math.} {\bf 138} (1993), 601--624}
\medskip\newref{GuanSpruckI}{Guan B., Spruck J., The existence of hypersurfaces of constant Gauss curvature with prescribed boundary, {\sl J. Differential Geom.} {\bf 62} (2002), no. 2, 259--287}
\newref{GuanSpruckII}{Guan B., Spruck J., Szapiel M., Hypersurfaces of constant curvature in Hyperbolic space I, {\sl J. Geom. Anal}}
\newref{GuanSpruckIII}{Guan B., Spruck J., Hypersurfaces of constant curvature in Hyperbolic space II, arXiv:0810.1781}
\newref{Gut}{Guti\'errez C., {\sl The Monge-Amp\`ere equation}, Progress in Nonlinear Differential Equations and Their Applications, {\bf 44}, Birkh\"a user, Boston, (2001)}
\newref{Jorgens}{J\"orgens K., \"Uber die L\"osungen der Differentialgleichung $rt-s^2=1$, {\sl Math. Ann.} {\bf 127} (1954), 130--134}
\newref{LabI}{Labourie F., Un lemme de Morse pour les surfaces convexes (French), {\sl Invent. Math.} {\bf 141} (2000), no. 2, 239--297}
\newref{LoftI}{Loftin J. C., Affine spheres and convex $\Bbb{RP}\sp n$-manifolds, {\sl Amer. J. Math.} {\bf 123} (2001), no. 2, 255--274}
\newref{LoftII}{Loftin J. C., Riemannian metrics on locally projectively flat manifolds, {\sl Amer. J. Math.} {\bf 124} (2002), no. 3, 595--609}
\newref{Pog}{Pogorelov A. V., On the improper affine hyperspheres, {\sl Geom. Dedicata} {\bf 1} (1972), 33--46}
\newref{PogB}{Pogorelov A. V., Extrinsic geometry of convex surfaces, Translated from the Russian by Israel Program for Scientific}
\newref{RosSpruck}{Rosenberg H., Spruck J., On the existence of convex hypersurfaces of constant Gauss curvature in hyperbolic space, {\sl J. Differential Geom.} {\bf 40} (1994), no. 2, 379--409}
\newref{SmiPPG}{Smith G., The Plateau problem for general curvature functions, arXiv:1008.3545}
\newref{SmiPPH}{Smith G., Compactness results for immersions of prescribed Gaussian curvature II - geometric aspects, arXiv:1002.2982}
\newref{SmiNLD}{Smith G., The Non-Linear Dirichlet Problem in Hadamard Manifolds, arXiv:0908.3590}
\newref{SmiFCS}{Smith G., Moduli of Flat Conformal Structures of Hyperbolic Type, arXiv:0804.0744}
\newref{SmiSLC}{Smith G., Special Lagrangian curvature, arXiv:math/0506230}
\newref{Spruck}{Spruck J., Fully nonlinear elliptic equations and applications to geometry, {\sl Proceedings of the International Congress of Mathematicians}, Vol. 1, 2 (Z\"urich, 1994), 1145--1152, Birkh\"auser, Basel, (1995)}
\newref{TrudWang}{Trudinger N. S., Wang X., On locally convex hypersurfaces with boundary, {\sl J. Reine Angew. Math.} {\bf 551} (2002), 11--32}

\ifmakebiblio%
    {\edef\fileoutput{\write\biblio{\egroup}}%
    \fileoutput}%
\fi%

%
%
%
\document
\myfontdefault
\global\chapno=1
\global\showpagenumflag=1
\def\Pagetitle{}
\null
\vfill
\def\centre{\rightskip=0pt plus 1fil \leftskip=0pt plus 1fil \spaceskip=.3333em \xspaceskip=.5em \parfillskip=0em \parindent=0em}%
\def\textmonth#1{\ifcase#1\or January\or Febuary\or March\or April\or May\or June\or July\or August\or September\or October\or November\or December\fi}
\font\abstracttitlefont=cmr10 at 14pt
{\abstracttitlefont\centre Compactness results for immersions\break of prescribed Gaussian curvature I - analytic aspects\par}
\bigskip
{\centre Graham Smith\par}
\bigskip
{\centre 9 December 2009\par}
\bigskip
{\centre Departament de Matem\`atiques,\par
Facultat de Ci\`encies, Edifici C,\par
Universitat Aut\`onoma de Barcelona,\par
08193 Bellaterra,\par
Barcelona,\par
SPAIN\par}
\bigskip
{\centre\it Dedicated to Hennie Smith\par
\medskip
22/12/1914-08/12/2009\par
\medskip
With affection.\par}
\bigskip
\noindent{\emph Abstract:\ }We extend recent results of Guan and Spruck, proving existence results for constant Gaussian curvature hypersurfaces in Hadamard manifolds.
%
%
\par
\vfill
\nextoddpage
\global\pageno=1
\def\Pagetitle{\sl Compactness Results I}
\newhead{Introduction}
\noindent Let $M:=M^{n+1}$ be an $(n+1)$-dimensional Riemannian manifold. An immersed hypersurface in $M$ is a pair $(\Sigma,\partial\Sigma):=((S,\partial S),i)$ where $(S,\partial S)$ is a compact, $n$-dimensional manifold with boundary and $i:S\rightarrow M$ is an immersion (that is, a smooth mapping whose derivative is everywhere injective). Throughout the sequel we abuse notation and denote $(S,\partial S)$ also by $(\Sigma,\partial\Sigma)$. We recall that the shape operator of the immersion is defined at each point by taking the covariant derivative in $M$ of the unit normal vector field over $\Sigma$ at that point, and that the {\bf Gaussian curvature} (also called the extrinsic curvature) is then defined as a function over $S$ to be equal to the determinant of the shape operator at each point.
\medskip
\noindent Geometers have studied the concept of Gaussian curvature ever since Gauss first proved in \cite{Gauss} his famous {\sl Teorema Egregium} which states that the Gaussian curvature of a surface immersed in $\Bbb{R}^3$ only depends on its intrinsic geometry and not on the immersion, which explains, for example, why a flat sheet of paper cannot be smoothly wrapped round a portion of the sphere. In more recent times, the Gaussian curvature of a hypersurface has revealed itself as an interesting object of study also from the perspective of geometric analysis as a straightforward and archetypal case of a much larger class of problems, including those of affine geometry, mass transport, Calabi-Yau geometry and so on, all of whose underlying equations are of so-called Monge-Amp\`ere type.
\medskip
\noindent In studying hypersurfaces of constant curvature of any sort, the most natural problems to study are those of Plateau and Minkowski, which ask respectively for the existence of hypersurfaces of constant curvature with prescribed boundary, or without boundary but instead satisfying certain topological conditions. The study of these problems has enjoyed a rich development over the last century, with the application of a wide variety of different techniques, including, for example, polyhedral approximation, used by Pogorelov to solve the Minkowski problem for convex, immersed spheres of prescribed Gaussian curvature in Euclidean space (c.f. \cite{PogB}), and, more recently, the continuity method, as used by Caffarelli, Nirenberg \& Spruck (c.f. \cite{CaffNirSprV}) to solve the Plateau problem for locally strictly convex (LSC) hypersurfaces which are graphs over a given hyperplane in Euclidean space. The ideas of Caffarelli, Nirenberg and Spruck were further developed in one direction by Rosenberg \& Spruck (c.f. \cite{RosSpruck}) to prove the existence of LSC hypersurfaces of constant extrinsic curvature in hyperbolic space with prescribed asymptotic boundary in the sphere at infinity (which was in turn generalised by Guan and Spruck in \cite{GuanSpruckII} and \cite{GuanSpruckIII} to treat more general notions of curvature). Likewise they were developed in another direction by Guan \& Spruck in \cite{GuanSpruckO} to prove existence of LSC hypersurfaces of constant extrinsic curvature in Euclidean space with prescribed boundary in the unit sphere. This led Spruck to conjecture in \cite{Spruck} that any compact, codimension $2$, immersed submanifold in Euclidean space which is the boundary of an LSC, immersed hypersurface is also the boundary of an LSC, immersed hypersurface of constant Gaussian curvature, a conjecture which was confirmed simultaneously by Guan \& Spruck in \cite{GuanSpruckI} and Trudinger \& Wang in \cite{TrudWang} using in both cases a combination of Caffarelli, Nirenberg \& Spruck's continuity method alongside an elegant application of the Perron method.
\medskip
\noindent With the exception of \cite{RosSpruck}, the above results essentially concern submanifolds of $\Bbb{R}^{n+1}$ and mild generalisations of this setup, and since most of the techniques used above rely in some way or another on the geometry of Euclidean space, the problem in general ambiant manifolds has remained largely open. Nonetheless, in \cite{LabI}, Labourie showed how pseudo-holomorphic geometry may be applied in conjunction with a parametric version of the continuity method to solve the Plateau problem in the case where $M$ is a $3$-dimensional Hadamard manifold. However, since this approach relies on techniques of holomorphic function theory, it does not easily generalise to the higher dimensional case, which has therefore hitherto remained unsolved. It is to fill this gap that we present in this and our forthcoming work \cite{SmiPPH} an approach which allows us to solve the Plateau problem for hypersurfaces of constant (or prescribed) Gaussian curvature in general manifolds, thus generalising the results \cite{GuanSpruckI} and \cite{TrudWang} of Guan \& Spruck and Trudinger \& Wang on the one hand and the result \cite{LabI} of Labourie on the other. In the interest of simplicity, we henceforth restrict attention to {\bf Hadamard manifolds}, which, we recall, are, by definition, complete, simply connected manifolds of non-positive sectional curvature. We leave the enthusiastic reader to investigate the few extra technical conditions required to state and prove the results in general manifolds.
\medskip
\noindent In the current paper, we will essentially be concerned with the local problem of finding solutions under conditions that are typically only valid over small regions. We will be mostly interested in the analysis required to obtain a-priori estimates and compactness results. In the forthcoming work, geometric results will be developed which will allow us to apply the estimates obtained here also to the global problem, as we will briefly discuss towards the end of this introduction.
\medskip
\noindent Thus let $M$ be a Hadamard manifold. Let $(\Sigma_0,\partial\Sigma_0)$ be a smooth, convex, immersed hypersurface in $M$ with smooth boundary. Let $\msf{N}$ be the exterior, unit, normal vector field over $\Sigma_0$ and define $\Cal{E}:\Sigma_0\times]-\infty,0]\rightarrow M$ by:
$$
\Cal{E}(x,t) = \opExp(-t\msf{N}).
$$
\noindent We have chosen here an unusual sign convention which we prefer for technical reasons. We say that a $C^{0,1}$ hypersurface $(\Sigma,\partial\Sigma)$ is a {\bf graph below} $\Sigma_0$ if and only if there exists a $C^{0,1}$ function $f:\Sigma_0\rightarrow]-\infty,0]$ and a homeomorphism $\varphi:\Sigma_0\rightarrow\Sigma$ such that:
\medskip
\myitem{(i)} $f$ vanishes along $\partial\Sigma_0$ (i.e. $\partial\Sigma=\partial\Sigma_0$); and
\medskip
\myitem{(ii)} for all $p\in\Sigma_0$:
$$
\varphi(p) = \opExp_p(-f(p)\msf{N}(p)).
$$
\noindent Let $(\hat{\Sigma},\partial\hat{\Sigma})$ be a $C^{0,1}$, convex, immersed hypersurface in $M$ which is a graph below $\Sigma_0$.
\medskip
\noindent We denote by $C_0^\infty(\Sigma_0)$ the space of smooth functions over $\Sigma_0$ which vanish along the boundary, and we identify surfaces which are graphs below $\Sigma_0$ with functions in $C_0^\infty(\Sigma_0)$. Gaussian curvature defines an operator $K:C_0^\infty(\Sigma_0)\rightarrow C^\infty(\Sigma_0)$ such that, for all $f\in C_0^\infty(\Sigma_0)$ and for all $p\in\Sigma$, $K(f)(p)$ is the Gaussian curvature of the graph of $f$ at the point below $p$. When the graph of $f$ is convex, the linearisation $DK_f$ of $K$ at $f$  is a second order, elliptic, partial differential operator. In particular, this is the case for $DK_0$, the linearisation of $K$ at the zero function, and we say that $(\Sigma_0,\hat{\Sigma})$ is {\bf stable} if and only if, for all $\psi\in C_0^\infty(\Sigma_0)$, if $DK_0\cdot\psi\geqslant 0$, then $\psi<0$ over the interior of $\Sigma_0$.
\medskip
\noindent We say that $(\Sigma_0,\hat{\Sigma})$ is {\bf rigid} if and only if there exists no other smooth hypersurface $\Sigma$ lying between $\Sigma_0$ and $\hat{\Sigma}$ such that $K(\Sigma) = K(\Sigma_0)$.
\medskip
\noindent In general, stable and rigid pairs of surfaces are relatively easy to construct inside small regions. For example, if $\Sigma_0$ is a bounded portion of a hypersurface in hyperbolic space which lies at constant distance from a totally geodesic hypersurface, then $(\Sigma_0,\hat{\Sigma})$ is both stable and rigid for any choice of $\hat{\Sigma}$. We refer the reader to Section \headref{HeadLocalAndGlobalRigidity} for more details.
\medskip
\noindent We prove the following local result:
\proclaim{Theorem \nextprocno}
\noindent Choose $k>0$ and suppose that the Gaussian curvature of $\Sigma_0$ is less than $k$. Suppose, moreover, that for some $\epsilon>0$ the Gaussian curvature of $\hat{\Sigma}$ is no less than $k+\epsilon$ in the weak (Alexandrov) sense and that the second fundamental form of $\hat{\Sigma}$ is also no less than $\epsilon$ in the weak (Alexandrov) sense. If $(\Sigma_0,\hat{\Sigma})$ is stable and rigid, then there exists a smooth, convex, immersed hypersurface $\Sigma_k$ such that:
\medskip
\myitem{(i)} $\Sigma_k$ is a graph below $\Sigma_0$;
\medskip
\myitem{(ii)} $\Sigma_k$ lies between $\Sigma_0$ and $\hat{\Sigma}$ as a graph below $\Sigma_0$; and
\medskip
\myitem{(iii)} the Gaussian curvature of $\Sigma_k$ is constant and equal to $k$.
\endproclaim
\proclabel{TheoremExistence}
\remark This follows immediately from Lemma \procref{LemmaExistence}.
\medskip
\remark The weak (Alexandrov) notion of lower (and upper) bounds for curvature is defined in Section \headref{HeadMaximumPrinciples}.
\medskip
\remark The hypothesis that $M$ be a Hadamard manifold is only made for simplicity of presentation. The same result, with appropriate modifications, continues to hold in more general manifolds.
\medskip
\noindent When $M$ is a space form, the Perron method may be applied to solve the following more general boundary value problem: let $\Gamma=(\Gamma_1,...,\Gamma_n)$ be a disjoint collection of closed, smooth, embedded, $(n-1)$-dimensional submanifolds of $\Bbb{H}^{n+1}$. Applying the machinery developed by Guan and Spruck in \cite{GuanSpruckI} along with Lemma \procref{LemmaDirichletInHyperbolicSpace} (which constitutes the more precise version of Theorem \procref{TheoremExistence} when $M=\Bbb{H}^{n+1}$) in place of Theorem $1.1$ of \cite{Guan}, we immediately obtain:
\proclaim{Theorem \nextprocno}
\noindent Choose $k>0$. Suppose that there exists a $C^2$, LSC, immersed hypersurface $\Sigma\subseteq\Bbb{H}^{n+1}$ of Gaussian curvature no less than $k$ such that $\partial\Sigma=\Gamma$. Then there exists a smooth (up to the boundary), locally strictly convex, immersed hypersurface $M\subseteq\Bbb{H}^{n+1}$ with $\partial M=\Gamma$ of constant Gaussian curvature equal to $k$. Moreover, $M$ is homeomorphic to $\Sigma$.
\endproclaim
\proclabel{TheoremBVP}
\noindent The proof of Theorem \procref{TheoremExistence} follows the analysis of Caffarelli, Nirenberg and Spruck first laid out in \cite{CaffNirSprI} and first applied to constant curvature hypersurfaces by the same authors in \cite{CaffNirSprV}. Our current work uses two key developments which simplify the analysis. The first, which is merely a question of perspective, is to analyse the Gauss Curvature Equation intrinsically along the hypersurface as in Section \headref{HeadSecondOrderBoundsOverTheInterior}, and the second is the use of Sard's Lemma in Section \headref{HeadOneDimensionalFamiliesOfSolutions} to generate smooth families of hypersurfaces interpolating between the data and the desired solution, which simplifies the topological approach already suggested by the work \cite{Guan} of Guan. In Section \headref{HeadRelationsToExistingResults} we show how our techniques can be easily adapted to recover both the results \cite{Guan} of Guan and \cite{RosSpruck} of Rosenberg and Spruck.
\medskip
\noindent As discussed previously, our main aim is to obtain a global existence result which confirms the natural extension of Spruck's conjecture (c.f. \cite{Spruck}) to more general manifolds. As we shall see in our forthcoming work \cite{SmiPPH}, the most significant new obstacle is the geometric problem of developing compactness results in general manifolds for LSC immersions with prescribed boundary. Having solved this problem, we return to Theorem \procref{TheoremExistence}, this time removing $\Sigma_0$ whilst also allowing $\hat{\Sigma}$ to vary between the data of the problem and a setup which may readily be shown to be stable and rigid. Then, proceeding as before, we obtain the following result, which is a mild simplifcation of the main result of \cite{SmiPPH}:
\proclaim{Theorem \nextprocno}
\noindent Let $M$ be a Hadamard manifold. Let $(\hat{\Sigma},\partial\hat{\Sigma})$ be a locally strictly convex, immersed hypersurface in $M$ whose boundary intersects itself transversally. Let $k>0$ be such that the Gaussian curvature of $\hat{\Sigma}$ is everywhere strictly greater than $k$.
\medskip
\noindent Suppose that there exists a convex subset $K\subseteq M$ with smooth boundary and an open subset $\Omega\subseteq\partial K$ such that:
\medskip
\myitem{(i)} $\partial\Omega$ is smooth; and
\medskip
\myitem{(ii)} $(\hat{\Sigma},\partial\hat{\Sigma})$ is isotopic through locally strictly convex immersions to a finite covering of $(\Omega,\partial\Omega)$,
\medskip
\noindent then there exists a locally strictly convex, immersed hypersurface $(\Sigma,\partial\Sigma)$ in $M$ such that:
\medskip
\myitem{(i)} $\partial\Sigma=\partial\hat{\Sigma}$; and
\medskip
\myitem{(ii)} $\Sigma$ has constant Gaussian curvature equal to $k$.
\endproclaim
\remark Since the submission of these papers, we have shown (c.f. \cite{SmiPPG}) that any locally strictly convex immersion is isotopic through locally strictly convex immersions to such a covering of an open subset of the boundary of a convex set, and so this condition is in fact redundant. We have chosen nonetheless to retain it here in order to keep this and the forthcoming paper as self-contained as possible.
\medskip
\noindent This paper is structured as follows:
\medskip
\myitem{(a)} in Section $2$, we show how first order bounds arise as a consequence of convexity;
\medskip
\myitem{(b)} in Section $3$, we derive the Gauss curvature equation for a graph in a general Riemannian manifold;
\medskip
\myitem{(c)} in Section $4$, we introduce the concept of weak (Alexandrov) lower and upper bounds for curvature;
\medskip
\myitem{(d)} in Sections $5$ and $6$ we obtain a-priori second order bounds over the boundary and then over the whole hypersurface respectively. These bounds are then applied in Section $7$ to obtain the compactness result, Lemma \procref{LemmaKisProper};
\medskip
\myitem{(e)} in Section $8$, we use Sard's Lemma to obtain smooth (albeit possibly empty) one-dimensional families of hypersurfaces interpolating between the data and the solutions. These are used in conjunction with the concepts of stability, rigidity and local rigidity developed in Sections $9$ and $10$ to prove in Section $10$ the existence result, Lemma \procref{LemmaExistence}, which immediately yields Theorem \procref{TheoremExistence};
\medskip
\myitem{(f)} In Section $11$, we restrict attention to space forms, proving Lemma \procref{LemmaDirichletInHyperbolicSpace}, which, in conjunction with the machinery developed by Guan and Spruck in \cite{GuanSpruckI} immediately yields Theorem \procref{TheoremBVP};
\medskip
\myitem{(g)} In Section $12$, we show how minor adaptations of these techniques allow us to obtain both the results \cite{Guan} of Guan (Theorem \procref{TheoremGuan}) and \cite{RosSpruck} of Rosenberg and Spruck (Theorem \procref{TheoremRosSpruck}); and
\medskip
\myitem{(h)} In Appendix $A$, we prove the regularity of limiting hypersurfaces which are themselves strictly convex. This result may be found in the notes of Caffarelli \cite{CaffII}, but given their general public unavailability, we consider it preferable to provide our own proof here.
\medskip
\noindent This paper was written whilst the author was staying at the Mathematics Department of the University Autonoma de Barcelona, Bellaterra, Spain.
\newhead{First Order Control}
\noindent Let $M:=M^{n+1}$ be an $(n+1)$-dimensional Riemannian manifold. Let $(\Sigma_0,\partial\Sigma_0)$ be a convex, immersed hypersurface with boundary. Let $\msf{N}_0$ and $A_0$ denote the outward pointing unit normal and the second fundamental form respectively of $\Sigma_0$. We define $\Cal{E}:\Sigma_0\times]-\infty,0]\rightarrow M$ by:
$$
\Cal{E}(x,t) = \opExp(-t\msf{N}_0(x)).
$$
\remark The change of sign ensures that convex hypersurfaces correspond to graphs of convex functions.
\medskip
\noindent We will say that a $C^{0,1}$ hypersurface, $\Sigma$, is a {\bf graph below} $\Omega$ if and only if there exists a $C^{0,1}$ function $f:\overline{\Omega}\rightarrow]-\infty,0]$ and a homeomorphism $\varphi:\overline{\Omega}\rightarrow\Sigma$ such that:
\medskip
\myitem{(i)} $f$ vanishes along $\partial\Omega$ (i.e. $\partial\Sigma=\partial\Omega$); and
\medskip
\myitem{(ii)} for all $p\in\Omega$:
$$
\varphi(p) = \opExp_p(-f(p)\msf{N}_0(p)).
$$
\noindent We refer to $f$ as the graph function of $\Sigma$. In particular, since $f$ is Lipschitz, its graph is never vertical, even along the boundary. Consider the family of graphs over $\Omega$. We define the partial order ``$<$'' on this family such that if $\Sigma$ and $\Sigma'$ are two graphs over $\Omega$ and $f$ and $f'$ are their respective graph functions, then:
$$
\Sigma < \Sigma' \Leftrightarrow f(p) < f'(p)\text{ for all }p\in\Omega.
$$
\noindent Since $\partial\Omega$ is smooth, for all $p\in\partial\Omega$, the set of supporting hyperplanes in $TM$ to $\partial\Omega$ at $p$ is parametrised by $\Bbb{R}$. Supporting hyperplanes may be locally considered as graphs over $\Omega$, and we obtain an analogous partial order on this set which we also denote by $<$.
\medskip
\noindent Let $\hat{\Sigma}$ be a $C^{0,1}$ convex hypersurface which is a graph over $\Omega$. Let $(\Sigma_n)_{\ninn}$ be a sequence of convex graphs over $\Omega$ such that for all $\ninn$, $\Sigma_n>\hat{\Sigma}$. For all $n$, let $f_n$ be the graph function of $\Sigma_n$.
\proclaim{Lemma \nextprocno}
\noindent $(f_n)_\ninn$ is uniformly bounded in the $C^{0,1}$ sense.
\endproclaim
\proclabel{LemmaFirstOrderBounds}
\proof For all $n\in\Bbb{N}\munion\left\{\infty\right\}$, define $U_n$ by:
$$
U_n = \left\{ \opExp_p(-t\msf{N}_0(p))\text{ s.t. }p\in\overline{\Omega}\text{ and }0\leqslant t\leqslant f_n(p)\right\}.
$$
\noindent By compactness of the family of convex sets, after extraction of a subsequence, there exists $U_0$ towards which $(U_n)_\ninn$ converges in the Hausdorff sense. Moreover, the supporting hyperplanes of $U_0$ are transverse to the normal geodesics leaving $H$. Indeed, suppose the contrary and let $p_0\in\partial U_0$ be a point where the supporting hyperplane is not transverse to the normal geodesic leaving $\Sigma_0$. Taking limits $\partial U_0\geqslant\hat{\Sigma}$. Since the tangent to $\hat{\Sigma}$ along $\partial\hat{\Sigma}$ is not vertical, it follows that $p_0$ lies over an interior point of $\Sigma_0$. Let $(p_n)_\ninn\in(\partial U_n)_\ninn$ be a sequence converging to $p_0$. For all $n\in\Bbb{N}\munion\left\{0\right\}$, let $q_n\in\Sigma_0$ be the orthogonal projection of $p_n$ onto $\Sigma_0$ and let $\gamma_n$ be the geodesic segment joining $q_n$ to $p_n$. For all $n\in\Bbb{N}$, $\gamma_n\subseteq U_n$. Taking limits, $\gamma_0\subseteq U_0$. It follows that $\gamma_0$ is an interior tangent to $\partial U_0$ at $p_0$. Therefore, by convexity, $\gamma_0\subseteq\partial U_0$. In particular, $U_0$ has a vertical supporting tangent at $q_0$, which is absurd. By compactness, we deduce that the supporting tangent hyperplanes of $(\partial U_n)_\ninn$ are uniformly transverse to the foliation of normal geodesics leaving $\Sigma_0$, and the result follows.\qed
\goodbreak
\newhead{The Gauss Curvature Equation}
\noindent Let $M:=M^{n+1}$ be an $(n+1)$-dimensional Riemannian manifold. Let $(\Sigma_0,\partial\Sigma_0)\subseteq M$ be a convex immersed hypersurface with boundary. Let $\msf{N}_0$ and $A_0$ denote the outward pointing unit normal and the second fundamental form respectively of $\Sigma_0$. Using the exponential map, we identify an open subset of $M$ with $\Sigma_0\times]-\infty,0]$.
\headlabel{HeadGaussCurvature}
\medskip
\noindent We will prove:
\goodbreak
\proclaim{Proposition \nextprocno}
\noindent Let $f:\Sigma_0\rightarrow]-\infty,0]$ be a smooth function. The Gaussian curvature of the graph of $f$ is given by:
$$
K = \psi(x,f,\nabla f)^{-1}\opDet(\opHess(f) + \Psi(x,f,\nabla f))^{1/n},
$$
\noindent where:
\medskip
\myitem{(i)} $\psi=\psi(x,t,p)$ is a smooth, strictly positive function and, for all $R>0$ there exists $\epsilon>0$ such that if $\left|t\right|<\epsilon$ then $\psi(x,t,p)$ is convex in $p$ for $\|p\|\leqslant R$; and
\medskip
\myitem{(ii)} there exists a smooth function $\Psi_0$ such that:
$$
\Psi(x,f,\nabla f)_{ij} = A_{0,ij}+f_{;i}f_{;k}{{A_0}^k}_j+f_{;j}f_{;k}{{A_0}^k}_i+f\Psi_0(x,f,\nabla f).
$$
\noindent Moreover, the graph of $f$ is convex if and only if $\opHess(f) + \Psi(x,f,\nabla f)$ is positive definite.
\endproclaim
\proclabel{PropFormulaForGaussCurvature}
{\sl\noindent Example:\ }We view $\Bbb{H}^n$ as a totally geodesic, embedded hypersurface in $\Bbb{H}^{n+1}$. Let $g_0$ and $g$ be the metrics of $\Bbb{H}^n$ and $\Bbb{H}^{n+1}$ respectively. We consider the foliation of $\Bbb{H}^{n+1}$ by geodesics normal to $\Bbb{H}^n$. Exceptionally, we reparametrise geodesics in a non-uniform manner in order to make this parametrisation conformal which simplifies the calculation of the connexion $2$-form. Let $\alpha:]-\pi/2,\pi/2[\rightarrow\Bbb{R}$ be such that, for all $\theta$:
$$
\opCos(\theta)\opCosh(\alpha(\theta)) = 1.
$$
\noindent Let $\msf{N}$ be the unit, normal vector field over $\Bbb{H}^n$ in $\Bbb{H}^{n+1}$. We define $\Phi:\Bbb{H}^n\times]-\pi/2,\pi/2[\rightarrow\Bbb{H}^{n+1}$ by:
$$
\Phi(x,\theta) = \opExp(-\alpha(t)\msf{N}(x)).
$$
\noindent We easily obtain:
$$
\Phi^* g = \frac{1}{\opCos^2(\theta)}(g_0\oplus d\theta^2).
$$
\noindent If $\Omega$ denotes the connexion $2$-form of the Levi-Civita covariant derivative of $\Phi^*g$ with respect to that of the product metric, then, for all $X$, $Y$ tangent to $\Bbb{H}^n$:
$$\matrix
\Omega(X,Y) \hfill&= -\langle X,Y\rangle\opTan(\theta)\partial_\theta,\hfill\cr
\Omega(X,\partial_\theta) \hfill&= \opTan(\theta)X,\hfill\cr
\Omega(\partial_\theta,\partial_\theta) \hfill&= \opTan(\theta)\partial_\theta.\hfill\cr
\endmatrix$$
\noindent Thus, if $\Omega\subseteq\Bbb{H}^n$ is an open set, and if $f:\Omega\rightarrow]-\pi/2,\pi/2[$ is a smooth function, then the Gaussian curvature of the graph of $f$ is given by:
$$
K = \opCos(f)^3(1 + \|\nabla f\|^2)^{-(n+2)/2n}\opDet(f_{;ij} - \opTan(f)(f_{;j}f_{;j} + \delta_{ij}))^{1/n}.
$$
\noindent We will return to this formula in later examples.\qed
\medskip
\noindent Let $\nabla^0$ denote the Levi-Civita covariant derivative of the product metric on $\Sigma_0\times ]-\infty,0]$. Let $g$ denote the pull back of the metric over $M$ through the exponential map. Let $\opVol$ denote the volume form of $g$ and let $\nabla$ denote the Levi Civita covariant derivative of $g$. Trivially, $\nabla$ coincides with the pull back through the exponential map of the Levi-Civita covariant derivative of $M$.
\proclaim{Proposition \nextprocno}
\noindent Let $\Omega:=\nabla-\nabla^0$ be the connection $2$-form of $\nabla$ with respect to $\nabla^0$. There exists a smooth $2$-form $\Omega_0$ such that, if $X$ and $Y$ are tangent to $\Sigma$, then:
$$\matrix
\Omega_{(x,t)}(X,Y) \hfill&= A_0(X,Y)\partial_t + t\Omega_{0,(x,t)}(X,Y),\hfill\cr
\Omega_{(x,t)}(X,\partial_t) \hfill&= -A_0 X + t\Omega_{0,(x,t)}(X,\partial_t),\hfill\cr
\Omega_{(x,t)}(\partial_t,\partial_t) \hfill&= t\Omega_{0,(x,t)}(\partial_t,\partial_t).\hfill\cr
\endmatrix$$
\endproclaim
\proclabel{PropConnexionInGraph}
\proof When $t=0$, by definition of $A_0$:
$$\matrix
\nabla_X Y \hfill&= \nabla^0_X Y + \langle \nabla_X Y,\msf{N}_0\rangle\msf{N}_0\hfill\cr
&= \nabla^0_X Y - A_0(X,Y)\msf{N}_0.\hfill\cr
\endmatrix$$
\noindent Thus, since $\msf{N}_0=-\partial_t$, at $t_0$:
$$
\nabla_X Y = \nabla^0_X Y + A_0(X,Y)\partial_t.
$$
\noindent Likewise:
$$
\nabla_X \partial_t = -\nabla_X \msf{N}_0 = -A_0 X.
$$
\noindent Finally, since the vertical lines are geodesics:
$$
\nabla_{\partial_t}\partial_t =0.
$$
\noindent The result follows.\qed
\medskip
\noindent Define $\hat{f}:\Sigma_0\times ]-\infty,0]\rightarrow\Bbb{R}$ by:
$$
\hat{f}(x,t) = f(x) - t.
$$
\noindent The graph of $f$ is the level set $\hat{f}^{-1}(\left\{0\right\})$. Observe that $\nabla\hat{f}$ is parallel to the downwards pointing unit normal over the graph of $f$. Let $A_f$ denote the second fundamental form of this graph. For all $i$, we define the vector field $\hat{\partial}_i = (\partial_i,f_{;i})_{(x,f(x))}$. $(\hat{\partial}_1,...,\hat{\partial}_n)$ forms a basis of the tangent space of the graph of $f$.
\medskip
{\bf\noindent Proof of Proposition \procref{PropFormulaForGaussCurvature}:\ } By definition:
$$
K^n = \opDet(A_f(\hat{\partial}_i,\hat{\partial}_j))/\opDet(g(\hat{\partial}_i,\hat{\partial}_j)).
$$
\noindent However, since the graph of $f$ is the level set $\hat{f}^{-1}(0)$:
$$
A_f = \frac{1}{\|\nabla\hat{f}\|_g}(\opHess(\hat{f})).
$$
\noindent Moreover:
$$
\opHess(\hat{f}) = \opHess^0(\hat{f}) - d\hat{f}(\Omega) = \opHess(f) - d\hat{f}(\Omega).
$$
\noindent It follows that $K$ has the specified form with:
$$
\psi(x,f,\nabla f) = \|\nabla\hat{f}\|_g\opDet(g(\hat{\partial}_i,\hat{\partial}_j))^{1/n},
$$
\noindent and:
$$
\Psi(x,f,\nabla f) = -d\hat{f}(\Omega).
$$
\noindent When $t=0$:
$$
\psi(x,0,p) = (1+\|p\|^2)^{(n+2)/2n}.
$$
\noindent Thus, since the function $p\mapsto(1+\|p\|^2)^\alpha$ is locally uniformly strictly convex for $\alpha>1/2$, $(i)$ follows.
\medskip
\noindent Likewise, by Proposition \procref{PropConnexionInGraph}:
$$\matrix
\Psi(x,0,p)(\hat{\partial}_i,\hat{\partial}_j) \hfill&= d\hat{f}(\msf{N})A_0(\partial_i,\partial_j) + f^{;j}d\hat{f}(A_0\partial_i) + f^{;i}d\hat{f}(A_0\partial_j)\hfill\cr
&= A_{0,ij} + f_{;i}f_{;k}{{A_0}^k}_j + f_{;j}f_{;k}{{A_0}^k}_i.\hfill\cr
\endmatrix$$
\noindent $(ii)$ follows.
\medskip
\noindent Finally, the graph of $f$ is convex if and only if $A_f$ is positive definite, and this completes the proof.\qed
\goodbreak
\newhead{Interlude - Maximum Principles}
\noindent Let $M:=M^{n+1}$ be an $(n+1)$-dimensional Riemanian manifold.
\headlabel{HeadMaximumPrinciples}
\proclaim{Definition \nextprocno}
\noindent Let $\Sigma$ be a $C^{0,1}$ convex, immersed hypersurface in $M$. Choose $k>0$. For $P\in\Sigma$, we say that the Gaussian curvature of $\Sigma$ is at least (resp. at most) $k$ in the weak (Alexandrov) sense at $P$ if and only if there exists a smooth, convex, immersed hypersurface $\Sigma'$ such that:
\medskip
\myitem{(i)} $\Sigma'$ is an exterior (resp. interior) tangent to $\Sigma$ at $P$; and
\medskip
\myitem{(ii)} the Gaussian curvature of $\Sigma'$ at $P$ is equal to $k$.
\endproclaim
\noindent This notion is well adapted to the weak Geometric Maximum Principle:
\proclaim{Lemma \nextprocno, {\bf Weak Geometric Maximum Principle}}
\noindent Let $\Sigma_1$, $\Sigma_2$ be two $C^{0,1}$, convex, immersed hypersurfaces in $M$. Choose $P\in\Sigma_1$. If $\Sigma_2$ is an interior tangent to $\Sigma_1$ at $P$, then the Gaussian curvature of $\Sigma_2$ at $P$ is no less than the Gaussian curvature of $\Sigma_1$ at $P$ in the weak (Alexandrov) sense.
\endproclaim
\proof Let $\Sigma_1'$ be a smooth, convex hypersurface which is an exterior tangent to $\Sigma_1$ at $P$. Likewise, let $\Sigma_2'$ be a smooth convex hypersurface which is an interior tangent to $\Sigma_2$ at $P$. Let $A_1$ and $A_2$ be the respective second fundamental forms of $\Sigma_1'$ and $\Sigma_2'$ respectively. Since $\Sigma_2'$ is an interior tangent to $\Sigma_1'$ at $P$:
$$
A_2 \geqslant A_1.
$$
\noindent The result follows.\qed
\medskip
\remark This result is often used in conjunction with foliations by constant curvature hypersurfaces which then act as barriers. In the case where $M=\Bbb{H}^{n+1}$, if we identify $\Bbb{H}^{n+1}$ with the upper half space in $\Bbb{R}^{n+1}$, then we obtain families of constant curvature hypersurfaces by considering intersections of spheres in $\Bbb{R}^{n+1}$ with $\Bbb{H}^{n+1}$. If the centre of such a sphere lies on $\Bbb{R}^n$, then its intersection with $\Bbb{H}^{n+1}$ has zero curvature. Otherwise, if the sphere is not entirely contained in $\Bbb{H}^{n+1}$, then the intersection has curvature less than $1$, and if it is contained in $\Bbb{H}^{n+1}$, then the intersection has curvature greater than $1$.
\medskip
\noindent We also have the strong Geometric Maximum Principle:
\proclaim{Lemma \nextprocno, {\bf Strong Geometric Maximum Principle}}
\noindent Let $(\Sigma_1,\partial\Sigma_1)$ and $(\Sigma_2,\partial\Sigma_2)$ be smooth, convex, immersed hypersurfaces in $M$ of constant Gaussian curvature equal to $k$.
\medskip
\myitem{(i)} If  $P$ is an interior point of $\Sigma_1$, and if $\Sigma_2$ is an exterior tangent to $\Sigma_1$ at $P$, then $\Sigma_1=\Sigma_2$.
\medskip
\myitem{(ii)} Suppose in addition that $\partial\Sigma_1=\partial\Sigma_2$. If $P$ is a boundary point of $\Sigma_1$ and if $\Sigma_2$ is an exterior tangent to $\Sigma_1$ at $P$, then $\Sigma_1=\Sigma_2$.
\endproclaim
\proof $\Sigma_2$ is a graph below $\Sigma_1$ near $P$. Let $U$ be a neighbourhood of $P$ in $\Sigma_1$ over which
$\Sigma_2$ is a graph. Let $A$ be the shape operator of $\Sigma_1$ and let $f$ be the graph function of $\Sigma_2$. By Proposition \procref{PropFormulaForGaussCurvature}:
$$
\opDet(\opHess(f) + \Psi(x,f,\nabla f))^{1/n}= k\psi(x,f,\nabla f),
$$
\noindent for some $\Psi$ and $\psi$. However:
$$
\opDet(A) = k.
$$
\noindent Thus, by concavity of $\opDet^{1/n}$:
$$
\frac{k}{n}\opTr(A^{-1}(\opHess(f) + \Psi(x,f,\nabla f) - A)) \geqslant k(\psi(x,f,\nabla f) - 1).
$$
\noindent Moreover, by the proof of Proposition \procref{PropFormulaForGaussCurvature}:
$$
\psi(x,f,\nabla f) = (1 + \|\nabla f\|^2)^{n+2/2n} + f\psi_0(x,f,\nabla f),
$$
\noindent For some smooth function $\psi_0$. Thus:
$$
k(\psi(x,f,\nabla f) - 1) = c_1f + \langle b_1,\nabla f\rangle,
$$
\noindent for some smooth function $c_1$ and vector field $b_1$. Likewise, by Proposition \procref{PropFormulaForGaussCurvature}:
$$
\opTr(A^{-1}(\Psi(x,f,\nabla f) - A)) = c_2 f + \langle b_2,\nabla f \rangle,
$$
\noindent for some smooth function $c_2$ and vector field $b_2$. Thus:
$$
\opTr(A^{-1}\opHess(f)) + \langle b,\nabla f\rangle + cf \geqslant 0,
$$
\noindent for some smooth function $c$ and vector field $b$. Since $f\leqslant 0$ and $f(P)=0$, in both cases $(i)$ and $(ii)$, it follows by the strong maximum principle (Theorems $3.5$ and $3.6$ of \cite{GilbTrud}) that $f=0$ over a neighbourhood of $P$. The result now follows by unique continuation of constant Gaussian curvature hypersurfaces.\qed
\goodbreak
\newhead{Second Order Bounds Along the Boundary}
\noindent Let $M:=M^{n+1}$ be an $(n+1)$-dimensional Riemannian manifold. Let $(\Sigma_0,\partial\Sigma_0)\subseteq M$ be a smooth, strictly convex, immersed hypersurface. Using the exponential map, we identify a subset of $M$ with $\Sigma_0\times]-\infty,0]$. Let $\phi:M\rightarrow ]0,\infty[$ be a smooth, positive function. Let $(\hat{\Sigma},\partial\hat{\Sigma})$ be a $C^{0,1}$, convex, immersed hypersurface such that:
\headlabel{HeadSecondOrderBoundaryBounds}
\medskip
\myitem{(i)} $\hat{\Sigma}$ is a graph below $\Sigma_0$;
\medskip
\myitem{(ii)} $\partial\hat{\Sigma}=\partial\Sigma_0$; and
\medskip
\myitem{(iii)} for all $x\in\hat{\Sigma}$, the Gaussian curvature of $\hat{\Sigma}$ is greater than $\phi(x)+\epsilon$ in the weak (Alexandrov) sense, for some $\epsilon>0$.
\medskip
\noindent $\hat{\Sigma}$ serves as a lower barrier for our problem. Let $(\Sigma,\partial\Sigma)\subseteq M$ be a smooth, convex, immersed hypersurface such that:
\medskip
\myitem{(i)} $\hat{\Sigma}\leqslant \Sigma\leqslant\Sigma_0$;
\medskip
\myitem{(ii)} $\partial\Sigma=\partial\Sigma_0$; and
\medskip
\myitem{(iii)} for all $x\in\Sigma$, the Gaussian curvature of $\Sigma$ at $x$ is equal to $\phi(x)$.
\medskip
\noindent We aim to obtain bounds for the norm of the second fundamental form of $\Sigma$ along the boundary which only depend on the data. To this end, we denote by $\Cal{B}$ the family of constants which depend continuously on the data: $M$, $\Sigma_0$, $\hat{\Sigma}$, $\epsilon$, $\phi$ and the $C^1$ jet of $\Sigma$ (formally, $\Cal{B}$ is the set of continuous - or even locally bounded - functions over the space of data). When supplementary data, $D$ (such as, for example, a vector field) is added, we denote by $\Cal{B}(D)$ the family of constants which, in addition, also depend on $D$.
\medskip
\noindent We will prove:
\proclaim{Proposition \nextprocno}
\noindent Let $\Sigma$, $\Cal{B}$ be as described above. If $\Sigma_0$ is strictly convex, then there exists $K$ in $\Cal{B}$ such that, if $A$ is the second fundamental form of $\Sigma$, then, for all $P\in\partial\Sigma$:
$$
\|A(P)\| \leqslant K.
$$
\endproclaim
\proclabel{PropSecondOrderBoundaryBounds}
\remark The strict convexity of $\Sigma_0$ is only required in the last step of the proof, where it is used to obtain uniform strict lower bounds for the restriction of the second fundamental form to the tangent space of $\partial\Sigma_0$. In other cases, such as where $\Sigma_0$ is totally geodesic, for example, this may shown using other means (c.f. Section \headref{HeadHyperbolicCase}).
\medskip
\noindent Let $P\in\partial\Sigma_0$ be a point on the boundary. For the sake of later applications (c.f. \cite{SmiPPH}), we underline that $\hat{\Sigma}$ need only exist locally. We thus let $\hat{\Sigma}_P\subseteq M$ be a smooth, convex, immersed hypersurface such that:
\medskip
\myitem{(i)} $\hat{\Sigma}_P$ is a graph below $\Sigma_0$;
\medskip
\myitem{(ii)} $P\in\hat{\Sigma}_P$; and
\medskip
\myitem{(iii)} for all $x\in\hat{\Sigma}_P$, the Gaussian curvature of $\hat{\Sigma}_P$ at $x$ is greater than $\phi(x)+\epsilon$.
\medskip
\noindent Bearing in mind the results of Section \headref{HeadGaussCurvature}, we will consider $\Sigma$ and $\hat{\Sigma}_P$ as graphs near $P$ over a hypersurface whose second fundamental form vanishes at $P$. Thus, let $\Sigma_1\subseteq M$ be an immersed hypersurface in $M$ which is tangent to $\Sigma_0$ at $P$ and which is totally geodesic at $P$.
\medskip
\noindent Let $\Omega\subseteq\Sigma_1$ be an open set with $P\in\partial\Omega$ and $f_0:\Omega\rightarrow\Bbb{R}$ a function such that:
\medskip
\myitem{(i)}  $\Sigma_0$ is the graph of $f_0$ over $\Omega$; and
\medskip
\myitem{(ii)} $f_0(\partial\Omega)=\partial\Sigma_0$.
\medskip
\remark Observe that both $\hat{\Sigma}_P$ and $\Sigma_1$ are local objects, only defined near $P$, as opposed to $\hat{\Sigma}$, for example, which is a global object, sharing the same boundary as $\Sigma_0$.
\medskip
\noindent We observe in passing that, by convexity, after reducing $\Sigma_1$ if necessary, $f_0$ may be made to be positive. $\partial\Omega$ consists of two components: we denote by $\partial_b\Omega$ the subset of $\partial\Omega$ which lies below the boundary of $\Sigma_0$ and we denote by $\partial_i\Omega$ the subset of $\partial\Omega$ which lies below the interior of $\Sigma_0$.
\proclaim{Proposition \nextprocno}
\noindent Let $\Sigma$, $\Cal{B}$ be as described at the beginning of this section. Let $\Omega$ be as described above. For all $P\in\partial\Sigma$, there exists $\delta>0$ in $\Cal{B}(P)$ and a neighbourhood $U$ of $P$ in $\Sigma$ which is a graph over $B_\delta(P)\minter\Omega$.
\endproclaim
\proof The radius over which $\Sigma$ is a graph over $\Sigma_1$ is determined by the $C^1$ jet of $\Sigma$, which is among the data defining $\Cal{B}$. The result follows.\qed
\medskip
\noindent We thus replace $\Omega$ with $\Omega\minter B_\delta(P)$ and let $f,\hat{f}:\Omega\rightarrow\Bbb{R}$ be the functions of which $\Sigma$ and $\hat{\Sigma}_P$ respectively are the graphs below $\Sigma_1$.
\medskip
\noindent By Proposition \procref{PropFormulaForGaussCurvature}, there exist functions $\psi$ and $\Psi$ and a positive number $R>0$, which only depends on $M$, $\phi$ and $\Sigma_1$ such that:
$$
\opDet(\opHess(f) + \Psi(x,f,\nabla f))^{1/n} = \psi(x,f,\nabla f).
$$
\noindent Moreover:
\medskip
\myitem{(i)} $\opHess(f) + \Psi(x,f,\nabla f)$ is positive definite;
\medskip
\myitem{(ii)} $\Psi(x,t,p), (\partial_{p_k}\Psi)(x,t,p)=O(d(x,P)) + O(t)$ where $d(\cdot,P)$ is the distance in $M$ to $P$; and
\medskip
\myitem{(iii)} for $t$ sufficiently small, $p \mapsto\psi(x,t,p)$ is a convex function in $p$ for $\|p\|\leqslant R$.
\medskip
\noindent We define the matrix $B$ by:
$$
B = \frac{1}{n}\psi(x,f,\nabla f)(\opHess(f) + \Psi(x,f,\nabla f))^{-1}.
$$
\noindent We define the operator $\Cal{L}$ by:
$$
\Cal{L}g = B^{ij}g_{;ij} + B^{ij}(\partial_{p_k}\Psi)_{ij}g_{;k} - (\partial_{p_k}\psi)g_{;k}.
$$
\proclaim{Proposition \nextprocno}
\noindent Let $\Sigma$, $\Cal{B}$ be as described at the beginning of this section. For all $P\in\partial\Sigma$, there exists $\delta_1>0$ and $\epsilon_1>0$ in $\Cal{B}(P)$ such that for $d(x,P)<\delta_1$:
$$
\Cal{L}(f-\hat{f}) \leqslant -\epsilon_1(1+\sum_{i=1}^n B^{ii}).
$$
\endproclaim
\proclabel{PropFirstBarrierEstimate}
\remark This inequality lies at the heart of the Caffarelli, Nirenberg, Spruck technique. The aim is to build functions which are subharmonic with respect to $\Cal{L}$, the key observation being that the appropriate term with respect to which bounds should be obtained is the trace of the matrix defining the generalised laplacian $\Cal{L}$, in this case $\sum_{i=1}^nB^{ii}$. We encourage the interested reader to compare this proposition with the relation shown on Line $13$ of Page $376$ of \cite{CaffNirSprI}, where the function $f-\hat{f}$ here plays the role of the function $x_n$ in the construction of their barrier function $w$. In addition, a clearer view of the main elements of the proof may be obtained by observing the effect of setting the constant $\eta_2$ to be equal to $0$, amounting to not perturbing $\hat{f}$. Finally, observe how the proof depends on the concavity of the determinant function as well as the convexity of $\psi(x,t,p)$ with respect to $p$, which is a recurring theme whenever this technique is applied.
\medskip
\proof There exists $\eta_1>0$ in $\Cal{B}$ such that, near $p$:
$$
\opDet(\opHess(\hat{f}) + \Psi(x,\hat{f},\nabla\hat{f}))^{1/n} \geqslant \psi(x,\hat{f},\nabla\hat{f}) + 2\eta_1.
$$
\noindent Define $\delta:\Sigma_1\rightarrow\Bbb{R}$ by:
$$
\delta(x) = d_1(x,P)^2,
$$
\noindent where $d_1$ denotes the intrinsic distance in $\Sigma_1$. Near $P$:
$$
\opHess(\delta) \geqslant \opId.
$$
\noindent There exists $\eta_2>0$ in $\Cal{B}$ such that, if we define $\hat{g}$ by:
$$
\hat{g} = \hat{f} - \eta_2\delta,
$$
\noindent then, near $P$:
$$
\opDet(\opHess(\hat{g}) + \Psi(x,\hat{g},\nabla\hat{g}))^{1/n} \geqslant \psi(x,\hat{g},\nabla\hat{g}) + \eta_1.
$$
\noindent Since $\opDet^{1/n}$ is a concave function:
$$\matrix
\opDet(\opHess(\hat{g}) + \Psi(x,\hat{g},\nabla\hat{g}))^{1/n}
-\opDet(\opHess(f) + \Psi(x,f,\nabla f))^{1/n}\hfill\cr
\qquad\leqslant B^{ij}(\hat{g}_{;ij} + \Psi_{ij}(x,\hat{g},\nabla\hat{g}) - f_{;ij} - \Psi_{ij}(x,f,\nabla f))\hfill\cr
\qquad\leqslant B^{ij}(\hat{f}-f)_{;ij} - \eta_2\sum_{i=1}^nB^{ii} + B^{ij}(\Psi_{ij}(x,\hat{g},\nabla\hat{g}) - \Psi_{ij}(x,f,\nabla f)).\hfill\cr
\endmatrix$$
\noindent Bearing in mind that $\Psi(x,t,p)=O(d(x,P)) + O(t)$, near $P$:
$$
B^{ij}(f-\hat{f})_{;ij} \leqslant - \eta_1 - \frac{\eta_2}{2}\sum_{i=1}^nB^{ii} + \psi(x,f,\nabla f) - \psi(x,\hat{g},\nabla\hat{g}).
$$
\noindent However, sufficiently close to $p$:
$$
\psi(x,f,\nabla\hat{g}) - \psi(x,\hat{g},\nabla\hat{g}) \leqslant \eta_1/3.
$$
\noindent Moreover, by convexity of $\psi$:
$$
\psi(x,f,\nabla f) - \psi(x,f,\nabla\hat{g}) \leqslant (\partial_{p_k}\psi)(f_{;k} - \hat{f}_{;k} + \eta_2\delta_{;k}).
$$
\noindent Since $\delta_{;k}$ is continuous and vanishes at $P$, we conclude that, near $P$:
$$
B^{ij}(f-\hat{f})_{;ij} - (\partial_{p_k}\psi)(f_{;k} - \hat{f}_{;k}) \leqslant -\eta_1/3 - \eta_2/2\sum_{i=1}^nB^{11}.
$$
\noindent Bearing in mind that, for all $k$, $(\partial_{p_k}\Psi)(x,t,\xi)=O(d(x,P)) + O(t)$, the result follows.\qed
\medskip
\noindent Let $X$ be a vector field over $\Sigma_1$.
\proclaim{Proposition \nextprocno}
\noindent Let $\Sigma$, $\Cal{B}$ be as described at the beginning of this section. For all $P\in\partial\Sigma$, there exists $K$ in $\Cal{B}(P,X)$ such that near $P$:
$$
\left|\Cal{L}(Xf)\right| \leqslant K(1+\sum_{i=1}^nB^{ii}).
$$
\endproclaim
\remark We encourage the interested reader to compare this relation to $2.12$ of \cite{CaffNirSprI}.
\medskip
\proof Differentiating the Gaussian curvature equation yields, for all $k$:
$$
B^{ij}(f_{;ijk} + (\partial_{x_k}\Psi)_{ij} + (\partial_t\Psi)_{ij}f_{;k} + (\partial_{p_l}\Psi)_{ij}f_{;lk}) = (\partial_{x_k}\psi) + (\partial_t\psi)f_{;k} + (\partial_{p_l}\psi)f_{;lk}.
$$
\noindent However:
$$
f_{;lk}=f_{;kl}.
$$
\noindent Moreover:
$$
f_{;ijk} = f_{;kij} + {R^{\Sigma_1}_{jki}}^pf_{;p},
$$
\noindent where $R^{\Sigma_1}$ is the Riemann curvature tensor of $\Sigma_1$. There therefore exists $K_1$ in $\Cal{B}(P,\Sigma_1)$ such that:
$$
\left|B^{ij}(f_{;kij} + (\partial_{p_l}\Psi)_{ij}f_{;kl}) - (\partial_{p_l}\psi)f_{;kl}\right| \leqslant K_1(1+\sum_{i=1}^nB^{ii}).
$$
\noindent Moreover, bearing in mind the definition of $B$, we obtain:
$$\matrix
B^{ij}f_{;ki} \hfill&= B^{ij}((f_{;ki} + \Psi_{ki}) - \Psi_{ki})\hfill\cr
&= \psi(x,f,\nabla f) - B^{ij}\Psi_{ki}.\hfill\cr
\endmatrix$$
\noindent However:
$$\matrix
\Cal{L}(Xf) \hfill&= X^k(B^{ij}(f_{;kij} + (\partial_{p_l}\Psi)_{ij}f_{;kl}) - (\partial_{p_l}\psi)f_{;kl})\hfill\cr
&\qquad + f_{;k}(B^{ij}({X^k}_{;ij} + (\partial_{p_l}\Psi)_{ij}{X^k}_{;l}) - (\partial_{p_l}\psi){X^k}_{;l})\hfill\cr
&\qquad + 2B^{ij}(f_{;ki}{X^k}_{;j}).\hfill\cr
\endmatrix$$
\noindent The result follows by combining the above relations.\qed
\proclaim{Corollary \nextprocno}
\noindent Let $\Sigma$, $\Cal{B}$ be as described at the beginning of this section. For all $P\in\partial\Sigma$, there exists $K$ in $\Cal{B}(P,X)$ such that near $P$:
$$
\left|\Cal{L}X(f-f_0)\right| \leqslant K(1+\sum_{i=1}^nB^{ii}).
$$
\endproclaim
\proclabel{CorSecondBarrierEstimate}
\noindent We define $\delta:\Sigma_1\rightarrow]0,\infty[$ by:
$$
\delta(x) = d_1(x,P)^2,
$$
\noindent where $d_1(\cdot,P)$ denotes the distance in $\Sigma_1$ to $P$.
\proclaim{Proposition \nextprocno}
\noindent Let $\Sigma$, $\Cal{B}$ be as described at the beginning of this section. For all $P\in\partial\Sigma$, there exists $K$ in $\Cal{B}(P,X)$ such that near $p$:
$$
\left|\Cal{L}\delta\right| \leqslant K(1 + \sum_{i=1}^nB^{ii}).
$$
\endproclaim
\proclabel{PropThirdBarrierEstimate}
\proof Trivial.\qed
\medskip
{\bf\noindent Proof of Proposition \procref{PropSecondOrderBoundaryBounds}:\ } Let $P$, $\Sigma_1$ and $\Omega$ be as before. Let $X$ be a vector field over $\Omega$ which is tangent to $\partial_b\Omega$. By Propositions \procref{PropFirstBarrierEstimate} and \procref{PropThirdBarrierEstimate} and Corollary \procref{CorSecondBarrierEstimate}, there exists $\eta,K>0$ in $\Cal{B}(P,X)$ such that:
$$\matrix
\left|\Cal{L}X(f-f_0)\right|\hfill&\leqslant K(1+\sum_{i=1}^nB^{ii}),\hfill\cr
\left|\Cal{L}\delta\right|\hfill&\leqslant K(1+\sum_{i=1}^nB^{ii}),\hfill\cr
\Cal{L}(f-\hat{f})\hfill&\leqslant -\eta(1+\sum_{i=1}^nB^{ii}).\hfill\cr
\endmatrix$$
\noindent Moreover, we may assume that, throughout $\Omega$:
$$
\left|X(f-f_0)\right| \leqslant K.
$$
\noindent By definition of $X$, $X(f-f_0)$ vanishes along $\partial_b\Omega$. Since $\partial_i\Omega$ is bounded away from $P$, there therefore exists $A_+>0$ in $\Cal{B}(P,X)$ such that, over $\partial\Omega$:
$$
X(f-f_0) - A_+\delta \leqslant 0.
$$
\noindent There exists $B_+>0$ in $\Cal{B}(P,X)$ such that, throughout $\Omega$:
$$
\Cal{L}(X(f-f_0) - A_+\delta - B_+(f-\hat{f})) \geqslant 0.
$$
\noindent Moreover, since $f-\hat{f}\geqslant 0$, this function is also negative along $\partial\Omega$. Thus, by the Maximum Principle, throughout $\Omega$:
$$
X(f-f_0) \leqslant A_+\delta + B_+(f-\hat{f}).
$$
\noindent Likewise, there exists $A_-,B_->0$ in $\Cal{B}(P,X)$ such that:
$$
X(f-f_0) \geqslant -A_-\delta - B_-(f-\hat{f}).
$$
\noindent There therefore exists $K_1>0$ in $\Cal{B}(P,X)$ such that:
$$
\left|d(X(f-f_0))(P)\right| \leqslant K_1.
$$
\noindent Thus, increasing $K_1$ if necessary:
$$
\left|d(Xf)(P)\right| \leqslant K_1.
$$
\noindent Let $\msf{N}$ be the unit normal vector field along $\partial\Sigma$ pointing into $\Sigma$. We have shown that there exists $K_2>0$ in $\Cal{B}$ such that, for any vector field, $X$, tangent to $\partial\Sigma$:
$$
\|A(X,\msf{N})\| \leqslant K_2\|X\|.
$$
\noindent The restriction of $A$ to $\partial\Sigma$ is determined by the norm of the second fundamental form of $\partial\Sigma=\partial\Sigma_0$. There therefore exists $K_3>0$ in $\Cal{B}$ such that, if $X$ and $Y$ are vector fields tangent to $\partial\Sigma$, then:
$$
\|A(X,Y)\|\leqslant K_3\|X\|\|Y\|.
$$
\noindent Finally, since $\Sigma$ lies between $\Sigma_0$ and $\hat{\Sigma}$, both of which are strictly convex, there exists $\epsilon_1>0$ in $\Cal{B}$ such that, throughout $\partial\Sigma$:
$$
A|_{T\partial\Sigma} \geqslant \epsilon_1\opId.
$$
\noindent Since $\opDet(A)=\phi$, $A(\msf{N},\msf{N})$ may be estimated in terms of the other components of $A$, and there therefore exists $K_4>0$ in $\Cal{B}$ such that, throughout $\partial\Sigma$:
$$
\|A(\msf{N},\msf{N})\|\leqslant K_4.
$$
\noindent The result now follows.\qed
\goodbreak
\newhead{Second Order Bounds Over the Interior}
\noindent Let $M:=M^{n+1}$ be a Hadamard manifold. Let $\Omega\subseteq M$ be a relatively compact open subset. Let $(\Sigma,\partial\Sigma)\subseteq M^{n+1}$ be a smooth, convex hypersurface and suppose that $\Sigma\subseteq\Omega$. Let $\msf{N}$ and $A$ be the unit, exterior, normal vector and the shape operator of $\Sigma$ respectively. In this section, it will be convenient to use the logarithm of the extrinsic curvature. Let $\phi:M\rightarrow\Bbb{R}$ be a strictly positive smooth function. We prove global second order estimates given second order estimates along the boundary for the problem:
\headlabel{HeadSecondOrderBoundsOverTheInterior}
$$
\opLog(\opDet(A)) = \phi(x),
$$
\noindent We denote by $\|A|_{\partial\Sigma_0}\|$ the supremum over $\partial\Sigma_0$ of the norm of $A$. We will prove:
\proclaim{Proposition \nextprocno}
\noindent There exists $K>0$ in $\Cal{B}(\|A|_{\partial\Sigma_0}\|)$ such that:
$$
\|A\| \leqslant K.
$$
\endproclaim
\proclabel{PropSecondOrderBounds}
\noindent In the sequel, we raise and lower indices with respect to $A$. Thus:
$$
A^{ij}A_{jk} = {\delta^i}_k,
$$
\noindent where $\delta$ is the Kr\"onecker delta function.
\proclaim{Proposition \nextprocno}
\myitem{(i)} For all $p$:
$$
A^{ij}A_{ij;p} = \phi_{;p}.
$$
\myitem{(ii)} For all $p,q$:
$$
A^{ij}A_{ij;pq} = A^{im}A^{jn}A_{ij;p}A_{mn;q} + \phi_{;pq}.
$$
\endproclaim
\proof This follows by differentiating the equation $\opLog(\opDet(A))=\phi$.\qed
\medskip
\noindent We recall the commutation rules of covariant differentiation in a Riemannian manifold:
\proclaim{Lemma \nextprocno}
\noindent Let $R^\Sigma$ and $R^M$ be the Riemann curvature tensors of $\Sigma$ and $M$ respectively. Then:
\medskip
\myitem{(i)} For all $i,j,k$:
$$
A_{ij;k} = A_{kj;i} + R^M_{ki\nu j},
$$
\noindent where $\nu$ represents the direction normal to $\Sigma$; and
\medskip
\myitem{(ii)} For all $i,j,k,l$:
$$
A_{ij;kl} = A_{ij;lk} + {R^\Sigma_{kli}}^pA_{pj} + {R^\Sigma_{klj}}^pA_{pi}.
$$
\endproclaim
\proclabel{LemmaCommutationRelations}
\proclaim{Corollary \nextprocno}
\noindent For all $i,j,k$ and $l$:
$$
A_{ij;kl} = A_{kl;ij} + {R^M_{kj\nu i}}_{;l} + {R^M_{li\nu k}}_{;j} + {R^\Sigma_{jlk}}^pA_{pi} + {R^\Sigma_{jli}}^pA_{pk}.
$$
\endproclaim
\proclabel{CorSecondDerivativesOfA}
\proof
$$\matrix
A_{ji;kl} \hfill&=A_{ki;jl} + {R^M_{kj\nu i}}_{;l}\hfill\cr
&=A_{ik;lj} + {R^M_{kj\nu i}}_{;l} + {R^\Sigma_{jlk}}^pA_{pi} + {R^\Sigma_{jli}}^pA_{pk}\hfill\cr
&=A_{lk;ij} + {R^M_{kj\nu i}}_{;l} + {R^M_{li\nu k}}_{;j} + {R^\Sigma_{jlk}}^pA_{pi} + {R^\Sigma_{jli}}^pA_{pk}\hfill\cr
\endmatrix$$
\noindent The result follows.\qed
\medskip
\noindent Choose $P\in\Sigma$. Let $\lambda_1,...,\lambda_n$ be the eigenvalues of $A$ at $P$. Choose an orthonormal basis, $(e_1,...,e_n)$ of $T_P\Sigma$ with respect to which $A$ is diagonal such that $a:=\lambda_1=A_{11}$ is the highest eigenvalue of $A$ at $P$. We extend this to a frame in a neighbourhood of $P$ by parallel transport along geodesics. We likewise extend $a$ to a function defined in a neighbourhood of $P$ by:
$$
a = A(e_1,e_1).
$$
\noindent Viewing $\lambda_1$ also as a function defined near $P$, $\lambda_1\geqslant a$ and $\lambda_1=a$ at $P$.
\proclaim{Proposition \nextprocno}
\noindent For all $i$, at $P$:
$$
a_{;ii} = A_{11;ii}.
$$
\endproclaim
\proof Bearing in mind that $\nabla_{e_i}e_i=0$ at $P$:
$$\matrix
a_{;ii} \hfill&= D_{e_i}D_{e_i}a\hfill\cr
&= D_{e_i}D_{e_i}A(e_1,e_1)\hfill\cr
&= D_{e_i}(\nabla A)(e_1,e_1;e_i) - 2D_{e_i}A(\nabla_{e_i}e_1,e_1)\hfill\cr
&= (\nabla^2A)(e_1,e_1;e_i,e_i) - 2A(\nabla_{e_i}\nabla_{e_i}e_1,e_1).\hfill\cr
\endmatrix$$
\noindent Since $e_1$ is defined by parallel transport along geodesics emanating from $P$, for all $i$, $\nabla_{e_i}\nabla_{e_i}e_1=0$ at $P$, and the result follows.\qed
\medskip
\noindent We define the Laplacian $\Delta$ such that, for all functions $f$:
$$
\Delta f = A^{ij}f_{;ij}.
$$
\proclaim{Proposition \nextprocno}
\noindent There exists $K>0$, which only depends on $M$ and $\phi$ such that, if $a>1$, then:
$$
\Delta\opLog(a)(P) \geqslant - K(1 + \sum_{i=1}^n\frac{1}{\lambda_i}).
$$
\endproclaim
\proclabel{PropSuperharmonic}
\proof By Corollary \procref{CorSecondDerivativesOfA}:
$$\matrix
a_{;ii} \hfill&= A_{11;ii}\hfill\cr
&= A_{ii;11} + {R^M_{i1\nu 1}}_{;i} + {R^M_{i1\nu i}}_{;1} + {R^\Sigma_{1ii}}^pA_{p1} + {R^\Sigma_{1i1}}^pA_{pi}.\hfill\cr
\endmatrix$$
\noindent However, at $P$:
$$
\sum_{i=1}^n\frac{1}{\lambda_1\lambda_i}A_{ii;11} = \sum_{i,j=1}^n\frac{1}{\lambda_i\lambda_j\lambda_1}A_{ij;1}A_{ij;1} + \frac{1}{\lambda_1}\phi_{;11}.
$$
\noindent Thus, at $P$:
$$\matrix
\Delta\opLog(a) \hfill&\geqslant \frac{1}{\lambda_1}\phi_{;11} +
\sum_{i,j=1}^n\frac{1}{\lambda_i\lambda_j\lambda_1}A_{ij;1}A_{ij;1} - \sum_{i=1}^n\frac{1}{\lambda_1^n\lambda_i}A_{11;i}A_{11;i}\hfill\cr
&\qquad\qquad + \sum_{i=1}^n\frac{1}{\lambda_1\lambda_i}(R^M_{i1\nu 1;i}+R^M_{i1\nu i;1})\hfill\cr
&\qquad\qquad + \sum_{i,j=1}^n\frac{1}{\lambda_1\lambda_i}({R^\Sigma_{1ii}}^pA_{p1} + {R^\Sigma_{1i1}}^pA_{pi}).\hfill\cr
\endmatrix$$
\noindent We consider each contribution seperately. Since, for all $a,b\in\Bbb{R}$, $(a+b)^2\leqslant 2a^2+2b^2$, by Lemma \procref{LemmaCommutationRelations}, for all $i\geqslant 2$:
$$
A_{11;i}^2 = (A_{i1;1} + R^M_{i1\nu 1})^2 \leqslant 2A_{i1;1}^2 + 2(R^M_{i1\nu 1})^2
$$
\noindent Thus, bearing in mind that $\lambda_1\geqslant 1$, there exists $K_1$, which only depends on $M$ such that:
$$
\sum_{i,j=1}^n\frac{1}{\lambda_i\lambda_j\lambda_1}A_{ij;1}A_{ij;1} - \sum_{i=1}^n\frac{1}{\lambda_1^2\lambda_i}A_{11;i}A_{11;i}\geqslant -K_1\sum_{i=1}^n\frac{1}{\lambda_i}.
$$
\noindent For all $\xi$, $X$ and $Y$:
$$\matrix
\nabla^\Sigma\xi(Y;X) \hfill&= \nabla^M\xi(Y;X) - A(X,Y)\xi(N);\text{ and }\hfill\cr
X\xi(\msf{N}) \hfill&= \nabla^M\xi(\msf{N};X) + \xi(A X).\hfill\cr
\endmatrix$$
\noindent Thus:
$$\matrix
R^M_{i1\nu 1;i} \hfill&= (\nabla^MR^M)_{i1\nu 1;i} + \lambda_i(1 - \delta_{i1}) R^M_{1 \nu\nu 1} + \lambda_i R^M_{i1i1},\hfill\cr
R^M_{i1\nu i;1} \hfill&= (\nabla^MR^M)_{i1\nu i;1} - \lambda_1(1 - \delta_{i1}) R^M_{i \nu\nu i} - \lambda_1 R^M_{i1i1}.\hfill\cr
\endmatrix$$
\noindent Bearing in mind that $\lambda_1\geqslant 1$, there exists $K_3$, which only depends on $M$ such that:
$$
\sum_{i=1}^n\frac{1}{\lambda_1\lambda_i}(R^M_{i1\nu 1;i}+R^M_{i1\nu i;1})\geqslant -K_3(1 + \sum_{i=1}^n\frac{1}{\lambda_i}).
$$
\noindent Moreover:
$$
{R^\Sigma_{1ii}}^pA_{p1} + {R^\Sigma_{1i1}}^pA_{pi} = R^M_{1ii1}(\lambda_1-\lambda_i) + \lambda_1\lambda_i(\lambda_1 - \lambda_i).
$$
\noindent Bearing in mind that $\lambda_1\geqslant 1$ and that $\lambda_1\geqslant\lambda_i$ for all $i$, there exists $K_2$, which only depends on $M$ such that:
$$
\sum_{i,j=1}^n\frac{1}{\lambda_1\lambda_i}({R^\Sigma_{1ii}}^pA_{p1} + {R^\Sigma_{1i1}}^pA_{pi}) \geqslant -K_2(1 + \sum_{i=1}^n\frac{1}{\lambda_i}).
$$
\noindent Since $\nabla^\Sigma_{e_1}e_1=0$ at $P$:
$$\matrix
\nabla^M_{e_1}e_1 \hfill&= \nabla^\Sigma_{e_1}e_1 + \langle\nabla^M_{e_1}e_1,\msf{N}\rangle\msf{N}\hfill\cr
&= -A(e_1,e_1)\msf{N}\hfill\cr
&= -\lambda_1\msf{N}.\hfill\cr
\endmatrix$$
\noindent Thus:
$$\matrix
\phi_{;11} \hfill&= \partial_1\partial_1\phi\hfill\cr
&= \opHess^M(\phi)(e_1,e_1) - d\phi(\nabla^M_{e_1}e_1)\hfill\cr
&= \opHess^M(\phi)(e_1,e_1) - \lambda_1d\phi(\msf{N}).\hfill\cr
\endmatrix$$
\noindent Bearing in mind that $\lambda_1\geqslant 1$, there thus exists $K_3$, which only depends on $M$ and $\phi$ such that:
$$
\frac{1}{\lambda_1}\phi_{;11} \geqslant -K_3.
$$
\noindent The result now follows by combining the above relations.\qed
\medskip
\noindent We recall that a function $f$ is said to satisfy $\Delta f\geqslant g$ in the weak sense if and only if, for all $P\in\Sigma$, there exists a smooth function $\varphi$, defined near $P$ such that:
\medskip
\myitem{(i)} $f\geqslant\varphi$ near $P$;
\medskip
\myitem{(ii)} $f=\varphi$ at $P$; and
\medskip
\myitem{(iii)} $\Delta\varphi\geqslant g$ at $P$.
\proclaim{Corollary \nextprocno}
\noindent With the same $K$ as in Proposition \procref{PropSuperharmonic}, if $\lambda_1\geqslant 1$, then:
$$
\Delta\opLog(\lambda_1) \geqslant -K(1+\sum_{i=1}^n\frac{1}{\lambda_i}),
$$
\noindent in the weak sense.
\endproclaim
\proclabel{CorSuperharmonic}
\proof Near $P\in\Sigma$, $\lambda_1\geqslant a$ and $\lambda_1=a$ at $P$. Since $P\in\Sigma_0$ is arbitrary, and since $a$ is smooth at $P$, the result follows.\qed
\medskip
\noindent Choose $x_0\in M$. Define $\delta$ by:
$$
\delta = \frac{1}{2}d(x,x_0)^2.
$$
\proclaim{Proposition \nextprocno}
\noindent There exists $c$, which only depends on $M$, $\Omega$, $\phi$ and $x_0$ such that:
$$
\lambda_1\geqslant c\ \Rightarrow\ \Delta^\Sigma\delta\geqslant \frac{1}{2}(1+\sum_{i=1}^n\frac{1}{\lambda_i}).
$$
\endproclaim
\proof Since $M$ has non-positive curvature:
$$\matrix
&\opHess^M(\delta)\hfill&\geqslant \opId\hfill\cr
\Rightarrow\hfill&\opHess^\Sigma(\delta)\hfill&\geqslant \opId - d(x,x_0)\langle\msf{N},\nabla d\rangle A\hfill\cr
\Rightarrow\hfill&\Delta \delta\hfill&\geqslant \sum_{i=1}^n\frac{1}{\lambda_i} - nd(x,x_0).\hfill\cr
\endmatrix$$
\noindent By compactness of $\Omega$, there exists $K_1>0$ such that, throughout $\Omega$:
$$
e^\phi \leqslant K_1.
$$
\noindent Thus:
$$\matrix
&\lambda_1\lambda_n^{n-1}\hfill&\leqslant K_1\hfill\cr
\Rightarrow\hfill&\lambda_n\hfill&\leqslant (K_1\lambda_1^{-1})^{1/(n-1)}\hfill\cr
\Rightarrow\hfill&\frac{1}{\lambda_n}\hfill&\geqslant (\lambda_1/K_1)^{1/(n-1)}\hfill\cr
\Rightarrow\hfill&\sum_{i=1}^n\frac{1}{\lambda_i}\hfill&\geqslant (\lambda_1/K_1)^{1/(n-1)}.\hfill\cr
\endmatrix$$
\noindent There thus exists $c_1>0$ such that, for $\lambda_1\geqslant c_1$, and for $x\in\Omega$:
$$\matrix
&\sum_{i=1}^n\frac{1}{\lambda_i}\hfill&\geqslant 2n\ d(x,x_0) + 1\hfill\cr
\Rightarrow\hfill&\Delta^\Sigma \delta\hfill&\geqslant \frac{1}{2}(1+\sum_{i=1}^n\frac{1}{\lambda}_i).\hfill\cr
\endmatrix$$
\noindent The result now follows.\qed
\proclaim{Corollary \nextprocno}
\noindent There exists $\lambda>0$ and $c>0$, which only depend on $M$, $\Omega$, $\phi$ and $x_0$ such that:
$$
\lambda_1\geqslant c\ \Rightarrow\ \Delta(\opLog(a) + \lambda \delta) > 0,
$$
\noindent in the weak sense.
\endproclaim
\proclabel{CorMaximumPrinciple}
\noindent Interior bounds now follow by the maximum principle:
\medskip
{\bf\noindent Proof of Proposition \procref{PropSecondOrderBounds}:\ } Consider the function $\|A\|e^{\lambda\delta} = \lambda_1e^{\lambda\delta}$. If this function achieves its maximum along $\partial\Sigma$, then the result follows since $e^{\lambda\delta}$ is uniformly bounded above and below. Otherwise, it acheives its maximum in the interior of $\Sigma$, in which case, by Corollary \procref{CorMaximumPrinciple} and the Maximum Principle, at this point:
$$
\|A\| = \lambda_1 \leqslant c.
$$
\noindent The result follows.\qed
\goodbreak
\newhead{Compactness}
\noindent Let $M:=M^{n+1}$ be an $(n+1)$-dimensional Hadamard manifold. Let $(\Sigma_0,\partial\Sigma_0)\subseteq M^{n+1}$ be a smooth, strictly convex hypersurface. Let $\msf{N}_0$ and $A_0$ be the unit, exterior, normal vector field and the shape operator of $\Sigma_0$ respectively. Using the exponential map, we identify $\Sigma\times]-\infty,0]$ with a subset of $M$.
\headlabel{HeadCompactness}
\medskip
\noindent Let $\opConv\subseteq C^\infty(\Sigma_0,]-\infty,0])$ be the family of smooth, negative valued functions over $\Sigma_0$ which vanish along $\partial\Sigma_0$ and whose graphs are strictly convex. We define the Gauss Curvature operator $K:\opConv\rightarrow C^\infty(\Sigma_0)$ such that, for all $f$, $(Kf)(x)$ is the Gauss curvature of the graph of $f$ at the point $(x,f(x))$. The formula for $K$ is given by Proposition \procref{PropFormulaForGaussCurvature}.
\medskip
\noindent Let $\hat{f}\in\opConv$ be such that:
$$
\hat{f}\leqslant 0,\qquad K(\hat{f}) - \epsilon > K(0) > 0,
$$
\noindent for some $\epsilon>0$. Denote $\phi_0=K(0)$ and $\hat{\phi}=K(\hat{f})$. Denote by $\opConv(\hat{f})$ the set of all $f\in\opConv$ such that:
$$
\hat{f}\leqslant f\leqslant 0,\text{ and }\hat{\phi}-\epsilon\geqslant K(f)\geqslant \phi_0.
$$
\noindent We prove a slightly stronger version of the assertion that the restriction of $K$ to $\opConv(\hat{f})$ is a proper mapping:
\proclaim{Lemma \nextprocno}
\noindent Let $(f_n)_\ninn$ be a sequence in $\opConv(\hat{f})$. Suppose there exists $(\phi_n)_\ninn\in C^\infty(M)$ such that, for all $n$, and for all $x\in\Sigma_0$:
$$
(Kf_n)(x) = \phi_n(x,f_n(x)).
$$
\noindent If there exists $\phi_\infty\in C^\infty(M)$ to which $(\phi_n)_\ninn$ converges, then there exists $f_\infty\in\opConv(\hat{f})$ to which $(f_n)_\ninn$ subconverges.
\endproclaim
\proclabel{LemmaKisProper}
\proclaim{Corollary \nextprocno}
\noindent The restriction of $K$ to $\opConv(\hat{f})$ is a proper mapping.
\endproclaim
\proclabel{CorKisProper}
{\bf\noindent Proof of Lemma \procref{LemmaKisProper}:\ } By Lemma \procref{LemmaFirstOrderBounds} and Propositions \procref{PropSecondOrderBoundaryBounds} and \procref{PropSecondOrderBounds}, there exists $C_1>0$ in $\Cal{B}$ such that, for all $n$:
$$
\|f_n\|_{C^2} \leqslant C_1.
$$
\noindent By Proposition \procref{PropFormulaForGaussCurvature}:
$$
Kf = F(\opHess(f), \nabla f, f, x),
$$
\noindent where $F(M,p,t,x)$ is elliptic in the sense of \cite{CaffNirSprII} and is concave in $M$. It follows by Theorem $1$ of \cite{CaffNirSprII} that there exists $\alpha>0$ and $C_2>0$ in $\Cal{B}$ such that, for all $n$:
$$
\|f_n\|_{C^{2,\alpha}} \leqslant C_2.
$$
\noindent Thus, by the Schauder Estimates (see \cite{GilbTrud}), for all $k\in\Bbb{N}$, there exists $B_k>0$ such that, for all $n$:
$$
\|f_n\|_{C^k} \leqslant B_k.
$$
\noindent The result now follows by the Arzela-Ascoli Theorem.\qed
\goodbreak
\newhead{One Dimensional Families of Solutions}
\noindent Let $M:=M^{n+1}$ be an $(n+1)$-dimensional Hadamard manifold. Let $(\Sigma_0,\partial\Sigma_0)\subseteq M^{n+1}$ be a smooth, convex hypersurface. Let $\msf{N}_0$ and $A_0$ be the unit, exterior, normal vector field and the shape operator of $\Sigma_0$ respectively. Using the Exponential Map, we identify $\Sigma\times]-\infty,0]$ with a subset of $M$.
\headlabel{HeadOneDimensionalFamiliesOfSolutions}
\medskip
\noindent Let $\hat{f}$,$\phi_0$ and $\hat{\phi}$ be as in the previous section. Let $\gamma:[0,1]\rightarrow C^\infty(\Sigma_0)$ be a smooth family of smooth functions such that, for all $\tau$:
$$
\phi_0 + \epsilon < \gamma(\tau) < \hat{\phi} -\epsilon,
$$
\noindent for some $\epsilon>0$. As before, let $K:\opConv\rightarrow C^\infty(\Sigma_0)$ be the Gauss Curvature Operator. For all $\phi\in C^\infty(\Sigma_0)$, define $\Gamma_\phi\subseteq I\times\opConv(\hat{f})$ by:
$$
\Gamma_\phi = \left\{ (t,f) \text{ s.t. }K(f) = \gamma(t) + \phi\right\}.
$$
\noindent  Viewing $\opConv$ as a Banach manifold (strictly speaking, the intersection of an infinite nested family of Banach manifolds), we will prove:
\proclaim{Proposition \nextprocno}
\noindent There exists $(\phi_n)_\ninn\in C^\infty(\Sigma_0)$ which converges to $0$ such that, for all $n$:
\medskip
\myitem{(i)} $\Gamma_n:=\Gamma_{\phi_n}$ is a (possibly empty) smooth, embedded, $1$-dimensional submanifold of $I\times\opConv(\hat{f})$; and
\medskip
\myitem{(ii)} $\partial\Gamma_n$ lies inside $\left\{0,1\right\}\times\opConv(\hat{f})$.
\endproclaim
\proclabel{PropTransversality}
\noindent We first prove:
\proclaim{Proposition \nextprocno}
\myitem{(i)} For all $\phi$, $\Gamma_\phi$ is compact; and
\medskip
\myitem{(ii)} For any neighbourhood $\Omega$ of $\Gamma_0$ in $I\times\opConv(\hat{f})$, there exists a neighbourhood $U$ of $0$ in $C^\infty(\Sigma_0)$ such that if $\phi\in U$, then $\Gamma_\phi\subseteq\Omega$.
\endproclaim
\proclabel{PropCompactnessOfGamma}
\proof $(i)$. This assertion follows from Corollary \procref{CorKisProper}.
\medskip
\noindent $(ii)$. Suppose the contrary. Let $(\tau_n)_\ninn\in [0,1]$, $(\phi_n)_\ninn\in C^\infty(\Sigma_0)$ and $(f_n)_\ninn\in\opConv(\hat{f})$ be such that $(\tau_n)_\ninn$ converges to $\tau_\infty\in[0,1]$, $(\phi_n)_\ninn$ converges to $0$ and, for all $n$:
$$
(\tau_n,f_n)\notin\Omega.
$$
\noindent Suppose moreover that, for all $n$:
$$
K(f_n) = \gamma(\tau_n) + \phi_n.
$$
\noindent By Lemma \procref{LemmaKisProper}, $(f_n)_\ninn$ subconverges to $f_\infty\in\opConv(f_\infty,\hat{f})$ such that:
$$\matrix
&K(f_\infty) \hfill&= \gamma(\tau_\infty)\hfill\cr
\Rightarrow\hfill& (\tau_\infty,f_\infty)\hfill&\in \Gamma_0.\hfill\cr
\endmatrix$$
\noindent Thus, for sufficiently large $n$, $(\tau_n,f_n)\in\Omega$, which is absurd. The result follows.\qed
\medskip
\noindent We denote by $C^\infty_0(\Sigma_0)$ the set of smooth functions on $\Sigma_0$ which vanish along $\partial\Sigma_0$, and we identify this with the tangent space of $\opConv$ in the natural manner. We consider the derivative of $K$:
\proclaim{Proposition \nextprocno}
\noindent At every point of $\opConv$, $DK$ defines a uniformly elliptic operator from $C_0^\infty(\Sigma_0)$ to $C^\infty(\Sigma_0)$.
\endproclaim
\proof This follows by differentiating the formula for the Gauss Curvature Operator given by Proposition \procref{PropFormulaForGaussCurvature}.\qed
\medskip
\noindent $DK$ is therefore Fredholm. Since it is defined on the space of smooth functions over a compact manifold with boundary, which themselves vanish over the boundary, it is of index zero.
\medskip
{\bf\noindent Proof of Proposition \procref{PropTransversality}:\ } Define $\hat{K}:[0,1]\times\opConv(\hat{f})\times C^\infty(\Sigma_0)\rightarrow C^\infty(\Sigma_0)$ by:
$$
\hat{K}(\tau,f,\phi) = \gamma(\tau) - K(f) + \phi.
$$
\noindent By compactness, there exists a neighbourhood $\Omega$ of $\Gamma$ in $[0,1]\times\opConv(\hat{f})$, a subspace $E\subseteq C^\infty(\Sigma_0)$ of dimension $m<\infty$ and $\epsilon>0$ such that the restriction of $D\hat{K}$ to $\Omega\times B_\epsilon(0)\subseteq\Omega\times E$ is always surjective. This restriction is Fredholm of index $(m+1)$. Define $\hat{\Gamma}$ by:
$$
\hat{\Gamma} = \hat{K}^{-1}(\left\{0\right\}).
$$
\noindent By the Implicit Function Theorem for Banach manifolds, $\hat{\Gamma}$ is an $(m+1)$-dimensional smooth submanifold of $\Omega\times B_\epsilon(0)$. Let $\pi_3:[0,1]\times\opConv(\hat{f})\times B_\epsilon(0)\rightarrow E$ denote projection onto the third factor. By Sard's Lemma there exists a sequence $(\phi_n)_\ninn\in B_\epsilon(0)$ which tends to $0$ such that, for all $n$, $\phi_n$ is a regular value of the restriction of $\pi_3$ to $\hat{\Gamma}$. However, for all $n$:
$$
\Gamma_n := \Gamma_{\phi_n} = \hat{\Gamma}\minter\pi_3^{-1}(\phi_n).
$$
\noindent Moreover, since $\phi_n$ is a regular value of $\pi_3$, $\Gamma_n$ is a (possibly empty) smooth $1$-dimensional embedded submanifold. By Proposition \procref{PropCompactnessOfGamma}, for all $n$, $\Gamma_n$ is compact, and for sufficiently large $n$, $\Gamma_n$ lies entirely inside $[0,1]\times\Omega$. Therefore:
$$
\partial\Gamma_n \subseteq \partial (I\times\Omega) \subseteq (\left\{0,1\right\}\times\Omega)\munion([0,1]\times\partial\opConv(\hat{f})).
$$
\noindent It thus remains to show that $\partial\Gamma_n$ lies away from $[0,1]\times\partial\opConv(\hat{f})$. However, if $(\tau,f)\in\Gamma_n$, then:
$$
0\leqslant f\leqslant\hat{f},\qquad \phi_0 + \epsilon < K(f) < K(\hat{f})-\epsilon.
$$
\noindent Thus, by the geometric maximum principle, away from $\partial\Sigma_0$:
$$
0 < f < \hat{f},
$$
\noindent and by the geometric maximum principal along the boundary, a similar relation holds for the derivative of $f$ in the internal normal direction along $\partial\Sigma_0$. It follows that $\Gamma_n$ lies in the interior of $\opConv(\hat{f})$ and so:
$$
\partial\Gamma_n\subseteq\left\{0,1\right\}\times\Omega.
$$
\noindent This completes the proof.\qed
\goodbreak
\newhead{Rigidity and Local Rigidity}
\noindent Let $M:=M^{n+1}$ be an $(n+1)$-dimensional Hadamard manifold. Let $(\Sigma_0,\partial\Sigma_0)\subseteq M^{n+1}$ be a smooth, convex hypersurface. Let $K:\opConv\rightarrow C^\infty(\Sigma_0)$, $\hat{f}$, $\phi_0$ and $\hat{\phi}$ be as in Section \headref{HeadCompactness}. Let $C_0^\infty(\Sigma_0)$ be the set of smooth functions over $\Sigma_0$, which, as in the preceeding section, we identify with the tangent space of $\opConv$. In particular, for all $f\in\opConv$ we denote by $DK_f:C_0^\infty(\Sigma_0)\rightarrow C^\infty(\Sigma_0)$ the derivative of $K$ at $f$.
\headlabel{HeadLocalAndGlobalRigidity}
\proclaim{Definition \nextprocno}
\myitem{(i)} We say that $\phi\in C^\infty(\Sigma_0)$ is locally rigid over $\opConv(\hat{f})$ if and only if for all $f\in\opConv(\hat{f})$ such that $K(f)=\phi$, $DK_f$ is invertible (in other words, $\phi$ is a regular value of $K$).
\medskip
\myitem{(ii)} We say that $\phi\in C^\infty(\Sigma_0)$ is rigid over $\opConv(\hat{f})$ if and only if there exists at most one $f\in\opConv(\hat{f})$ such that $K(f) = \phi$.
\endproclaim
\proclabel{DefnLocalAndGlobalRigidity}
{\sl\noindent Example:\ }Let $\Bbb{H}^{n+1}$ be $(n+1)$-dimensional hyperbolic space. Let $H$ be a totally geodesic hypersurface. For $D>0$, let $H(D)$ be the equidistant hypersurface at a distance $D$ from $H$. $H(D)$ has constant Gaussian curvature equal to $\opTanh(D)$. Let $\Omega\subseteq H(D)$ be any bounded open subset with smooth boundary and consider the hypersurface $(\Sigma_0,\partial\Sigma_0)=(\Omega,\partial\Omega)$. Define $f_0=0$ and $\phi_0=Kf_0=\opTanh(D)$. By the strong Geometric Maximum Principle and the homogeneity of $\Bbb{H}^{n+1}$, we readily show that $\phi_0$ is rigid for any choice of $\hat{\Sigma}$. Moreover, by calculating the Jacobi operator of $H(D)$ (or by calculating the derivative of $K$ using the example in Section \headref{HeadGaussCurvature}), we likewise show that $\phi_0$ is locally rigid.\qed
\medskip
{\sl\noindent Example:\ }The above example is a special case of a more general construction. Let $M$ be a Riemannian manifold. Let $P\in M$ be a point, let $\msf{N}\in UM$ be a unit vector at $P$, let $A$ be a positive-definite symmetric $2$-form over $\msf{N}^\perp$ and let $k=\opDet(A)$. There is no algebraic obstruction to the construction of a hypersurface $\Sigma$ such that:
\medskip
\myitem{(i)} $P\in\Sigma$;
\medskip
\myitem{(ii)} $\msf{N}$ is normal to $\Sigma$ at $P$;
\medskip
\myitem{(iii)} the second fundamental form of $\Sigma$ at $P$ is equal to $A$; and
\medskip
\myitem{(iv)} if $\psi=\opDet(A)$ is the Gaussian curvature of $\Sigma$, then $\psi=k$ up to infinite order at $P$.
\medskip
\noindent Since $\psi=k$ up to infinite order at $P$, for $\epsilon>0$ small, there exists a smooth family $(\psi_t)_{t<\epsilon}$ of smooth functions such that:
\medskip
\myitem{(i)} $\psi_0=\psi$; and
\medskip
\myitem{(ii)} for all $t$, $\psi_t=k$ over the geodesic ball of radius $t$ about $P$.
\medskip
\noindent Suppose moreover that $M$ has negative sectional curvature bounded above by $-1$ and that $A=k\opId$ for $k<1$. In this case, the derivative of the Gauss Curvature Operator is invertible over a geodesic ball of small radius about $P$ (see \cite{LabI} for details in the $2$-dimensional case). We may therefore assume by the Inverse Function Theorem for Banach Manifolds that $\psi=k$ over a geodesic ball of small radius about $P$. Moreover, $\Sigma$ may be extended to a foliation $(\Sigma_t)_{t\in]-\epsilon,\epsilon[}$ of a neighbourhood of $P$ in $M$ by hypersurfaces of constant curvature equal to $k$. Now let $B\subseteq M$ be a geodesic ball in $M$ centred on $P$ which is covered by this foliation. Let $\Omega\subseteq\Sigma$ be an open set with smooth boundary contained in $B\minter\Sigma$. If $\Sigma'$ is any other hypersurface of constant Gaussian curvature equal to $k$ such that $\partial\Sigma'=\partial\Omega$, then, by the Geometric Maximum Principle, $\Sigma'$ is contained inside $B$, and, by the strong Geometric Maximum Principle, $\Sigma'$ coincides with a leaf of the foliation. It is therefore equal to $\Omega$, and we have thus shown that $\psi=k$ is both rigid and locally rigid over $\Omega$ for any choice of $\hat{\Sigma}$.\qed
\proclaim{Proposition \nextprocno}
\myitem{(i)} If $\phi$ is locally rigid, then $\phi'$ is also locally rigid for all $\phi'$ sufficiently close to $\phi$.
\medskip
\myitem{(ii)} If $\phi$ is rigid and locally rigid, then $\phi'$ is rigid for all $\phi'$ sufficiently close to $\phi$.
\endproclaim
\proclabel{PropPropertiesOfRigidity}
\proof $(i)$. Suppose the contrary. Let $(\phi_n)_\ninn\in C^\infty(\Sigma_0)$ be a sequence of non-locally rigid functions converging to $\phi$. Since $\phi$ is locally rigid, $DK$ is invertible at $f$ for all $f\in K^{-1}(\left\{\phi\right\})$. There therefore exists a neighbourhood $\Omega$ of $K^{-1}(\left\{\phi\right\})$ in $\opConv(\hat{f})$ such that $DK$ is invertible at $f$ for all $f\in\Omega$. However, by Corollary \procref{CorKisProper}, for all sufficiently large $n$:
$$
K^{-1}(\left\{\phi_n\right\})\subseteq\Omega.
$$
\noindent $\phi_n$ is therefore locally rigid for sufficiently large $n$, which is absurd, and the assertion follows.
\medskip
\noindent $(ii)$. Suppose the contrary. There exists a sequence $(\phi'_n)_\ninn$ which converges to $\phi$ such that $\phi_n'$ is not globally rigid. Thus, for all $n$, there exists $f_{1,n}\neq f_{2,n}\in\opConv(\hat{f})$ such that:
$$
Kf_{1,n} = Kf_{2,n} = \phi_n'.
$$
\noindent By Corollary \procref{CorKisProper}, there exist $f_{1},f_{2}\in\opConv(\hat{f})$ to which $(f_{1,n})_\ninn$ and $(f_{2,n})_\ninn$ respectively converge. In particular:
$$
Kf_{1} = Kf_{2} = \phi.
$$
\noindent Since $\phi$ is rigid, it follows that:
$$
f_{1} = f_{2} = f.
$$
\noindent Since $\phi$ is locally rigid, $DK$ is invertible at $f$ and thus $K$ is locally invertible over a neighbourhood of $f$. In particular, for sufficiently large $n$:
$$
f_{1,n} = f_{2,n}.
$$
\noindent This is absurd, and the result follows.\qed
\goodbreak
\newhead{Stability and Existence}
\noindent Let $M:=M^{n+1}$ be an $(n+1)$-dimensional Hadamard manifold. Let $(\Sigma_0,\partial\Sigma_0)\subseteq M^{n+1}$ be a smooth, convex hypersurface. Let $(\hat{\Sigma},\partial\hat{\Sigma})\subseteq M^{n+1}$ be another smooth, convex hypersurface which is a graph below $\Sigma_0$. Let $\hat{f}\in C_0^\infty(\Sigma_0)$ be the function of which $\hat{\Sigma}$ is a graph. As in Section \headref{HeadCompactness}, we denote $\hat{\phi} = K(\hat{f})$ and $\phi_0=K(0)$, and we denote by $\opConv(\Sigma_0,\hat{\Sigma}):=\opConv(\hat{f})$ the set of all smooth functions in $C_0^\infty(\Sigma_0)$ such that:
$$
\hat{f} \leqslant f\leqslant 0,\ \text{and}\ \hat{\phi}-\epsilon\geqslant K(f)\geqslant \phi_0.
$$
\noindent We identify every function in $\opConv(\Sigma_0,\hat{\Sigma})$ with its graph.
\proclaim{Definition \nextprocno}
\myitem{(i)} We say that $(\Sigma_0,\hat{\Sigma})$ is stable if and only if for all $\psi\in C_0^\infty(\Sigma_0)$, if $DK_0\psi>0$, then $\psi<0$ over the interior of $\Sigma_0$.
\medskip
\myitem{(ii)} We say that $(\Sigma_0,\hat{\Sigma})$ is rigid if and only if the only hypersurface $(\Sigma,\partial\Sigma)\in\opConv(\Sigma_0,\hat{\Sigma})$ such that $K(\hat{\Sigma}) = K(\Sigma_0)$ is $\Sigma_0$ itself.
\endproclaim
\remark In other words, $(\Sigma_0,\hat{\Sigma})$ is rigid if and only if $\phi_0:=K(0)\in C^\infty(\Sigma_0)$ is rigid over $\opConv(\hat{f})$.
\medskip
\remark Observe that if $(\Sigma_0,\hat{\Sigma})$ is both rigid and stable, then $\phi_0$ is also locally rigid over $\opConv(\hat{f})$.
\medskip
{\noindent\sl Example:\ }Let $\msf{N}_0$ and $A_0$ be respectively the outward pointing, unit, normal vector field over $\Sigma_0$ and its corresponding shape operator. Let $JK$ be the Jacobi operator of $\Sigma_0$. $JK$ measures the first order variation of the Gaussian curvature upon first order, normal perturbations of $\Sigma_0$ and is given by:
$$
JK\phi = \opTr(A_0^{-1}W - A_0)\phi - \opTr(A_0^{-1}\opHess(\phi)).
$$
\noindent where the mapping $W$ is given by:
$$
\langle W(X),Y\rangle = \langle R_{X\msf{N}_0}Y,\msf{N}_0\rangle,
$$
\noindent and where $R$ is the Riemann curvature tensor of $M$. It follows that if the sectional curvature of $M$ is bounded above by $-\epsilon^2$ and if the principal curvatures of $\Sigma_0$ are bounded below by $\epsilon$, then:
$$
JK\phi = h\phi - \opTr(A_0^{-1}\opHess(\phi)),
$$
\noindent for some non-negative function $h$. Since $DK$ is conjugate to $JK$, it follows from the maximum principal that $(\Sigma_0,\hat{\Sigma})$ is stable.\qed
\proclaim{Lemma \nextprocno}
\noindent If $(\Sigma_0,\hat{\Sigma})$ is stable and rigid, then, for all $\phi$ such that:
$$
\phi_0 \leqslant \phi \leqslant \hat{\phi}-\epsilon,
$$
\noindent for some $\epsilon>0$, there exists a smooth, convex, immersed hypersurface $\Sigma_\phi$ such that:
\medskip
\myitem{(i)} $\hat{\Sigma}\leqslant\Sigma_\phi\leqslant\Sigma_0$, and
\medskip
\myitem{(ii)} the Gaussian curvature of $\Sigma_\phi$ at the point $p$ is equal to $\phi(p)$.
\endproclaim
\proclabel{LemmaExistence}
\proof Assume first that:
$$
\phi_0 + \epsilon < \phi < \hat{\phi}-\epsilon.
$$
\noindent By stability, reducing $\epsilon$ is necessary, there exists $f_0\in\opConv(\hat{f})$ such that:
$$
\phi_0' := K(f_0) > \phi_0 + \epsilon.
$$
\noindent By Proposition \procref{PropPropertiesOfRigidity}, we may assume moroever that $\phi_0'$ is both rigid and locally rigid over $\opConv(\hat{f})$. Let $\gamma:[0,1]\rightarrow C^\infty(\Sigma_0)$ be a smooth family of smooth functions such that:
\medskip
\myitem{(i)} $\gamma(0)=\phi_0'$, $\gamma(1)=\phi$, and
\medskip
\myitem{(ii)} for all $t\in[0,1]$:
$$
\phi_0 + \epsilon < \gamma(t) < \hat{\phi} - \epsilon.
$$
\noindent By Proposition \procref{PropTransversality}, there exists $(\phi_n)_\ninn\in C^\infty(\Sigma_0)$ which converges to $0$ such that, for all $n$, $\Gamma_n:=\Gamma_{\phi_n}$ is a (possibly empty) smooth, $1$-dimensional embedded submanifold of $[0,1]\times\opConv(\hat{f})$. Moreover, for all $n$, $\Gamma_n$ is compact, and:
$$
\partial\Gamma_n \subseteq \left\{0,1\right\}\times\opConv(\hat{f}).
$$
\noindent By Proposition \procref{PropPropertiesOfRigidity}, we may assume that, for all $n$, $\phi_0'+\phi_n$ is both rigid and locally rigid. Likewise, since $\phi_0'$ is locally rigid, we may assume that, for all $n$, there exists $f_n\in\opConv(\hat{f})$ such that:
$$
(0,f_n) \in \Gamma_n.
$$
\noindent $\Gamma_n$ is therefore non-empty for all $n$. Let $\Gamma_n^0$ be the connected component of $\Gamma_n$ containing $(0,f_n)$. Since it is compact, it is either an embedded, compact interval or an embedded, closed loop. We claim that $\Gamma_n^0$ is not a closed loop. Indeed, by local rigidity, $DK$ is invertible at $(0,f_n)$. Consequently, if $\pi_1:[0,1]\times\opConv(\hat{f})\rightarrow[0,1]$ is the projection onto the first factor, the restriction of $D\pi_1$ to $T\Gamma^0_n$ is invertible at $f_n$. Since $0=\pi_1(f_n)$ is an end point of $[0,1]$, this would imply that $(0,f_n)$ is also an end point of $\Gamma_n$. This is absurd and the assertion follows.
\medskip
\noindent For all $n$, let $g_n$ by the other end of $\Gamma_n^0$. Since $(\phi_0'+\phi_n)$ is globally rigid:
$$
g_n \in \left\{1\right\}\times\opConv(\hat{f}).
$$
\noindent In other words:
$$
K(g_n) = \phi + \phi_n.
$$
\noindent By Corollary \procref{CorKisProper}, there exists $g_0\in\opConv(\hat{f})$ to which $(g_n)_\ninn$ subconverges. In particular:
$$
K(g_0) = \phi.
$$
\noindent This proves existence in the case where $\phi_0+\epsilon<\phi<\hat{\phi}-\epsilon$. The general case follows by taking limits.\qed
\goodbreak
\newhead{Space Forms and the Local Geodesic Condition}
\noindent Let $M:=M^{n+1}$ be an $(n+1)$-dimensional Riemannian manifold. Let $K\subseteq M$ be a convex set with non-trivial interior. For any $P\in\partial K$, we say that $K$ satisfies the {\bf local geodesic condition} at $P$ if and only if there exists an open geodesic segment $\Gamma$ such that:
\headlabel{HeadHyperbolicCase}
\medskip
\myitem{(i)} $P\in\Gamma$; and
\medskip
\myitem{(ii)} $\Gamma\subseteq K$.
\medskip
\noindent Observe that since $K$ is convex, the second condition implies in particular that $\Gamma\subseteq\partial K$.
\medskip
\noindent We henceforth restrict attention to the case where $M$ is a space-form of non-positive curvature. In other words, up to rescaling, $M$ is isometric to either $(n+1)$-dimensional Euclidean or Hyperbolic space. We obtain the following global consequence of the local geodesic condition (c.f. \cite{SmiFCS}):
\proclaim{Lemma \nextprocno}
\noindent Let $K$ be a bounded, convex set with non-trivial interior, let $X\subseteq\partial K$ be a closed subset and let $Y\subseteq\partial K$ be the set of all points in $\partial K$ satisfying the local geodesic condition. If $X\munion Y$ is closed, then $Y$ lies in the convex hull of $X$.
\endproclaim
\proclabel{LemmaClosedMeansInConvexHull}
\proof Suppose that $Y$ is not contained in the convex hull of $X$. Then there exists a point $Q\in Y$ and a supporting hyperplane $H$ to $K$ at $Q$ such that $H\minter X$ is empty. Denote $K'=K\minter H$. Let $Y'\subseteq\partial K'$ be the set of points satisfying the local geodesic condition. In particular, $Q\in Y'$. Since $X\munion Y$ is closed, so is $Y'$. We now show that $Y'$ is unbounded. Indeed, suppose the contrary. Choose any $P\in M$ and let $d_P$ be the distance to $P$ in $M$. Since $Y'$ is closed and bounded, it is compact, and so there exists a point $Q'\in Y'$ maximising $d_P$. Let $\Gamma$ be the open geodesic segment in $Y'$ passing through $Q'$. Trivially, $\Gamma\subseteq Y'$. However, the restriction of $d_P$ to $\Gamma$ is convex, and so it cannot have a local maximum at $Q'$. This is absurd, and the assertion follows. However, since $K$ is bounded, so is $Y'$. This is absurd, and the result follows.\qed
\medskip
\noindent In the current context, regularity follows from the following result:
\proclaim{Proposition \nextprocno}
\noindent Suppose that $M$ is a space form of non-positive curvature. Choose $k>0$ and let $(K_n)_\ninn\subseteq M$ be a sequence of convex subsets of $M$ with smooth boundary such that, for all $n$, the Gaussian curvature of $\partial K_n$ is equal to $k$. Suppose that $(K_n)_\ninn$ converges to $K_0\subseteq M$ and that $K_0$ has non-empty interior. Then the set of points in $\partial K_n$ satisfying the local geodesic condition is closed.
\endproclaim
\proclabel{PropNodalSetIsClosed}
\proof We show that the complement is open. Indeed, let $Q\in\partial K_n$ be a point not satisfying the local geodesic condition. Then there exists a hyperplane $H$, a bounded, open, convex subset $U$ of $H$ and an open subset $V$ of $\partial K_0$ such that:
\medskip
\myitem{(i)} $Q\in V$;
\medskip
\myitem{(ii)} $Q$ lies at non-zero distance from $H$; and
\medskip
\myitem{(iii)} $V$ is a graph over $U$ with $\partial V=\partial U$.
\medskip
\noindent It follows from \cite{CaffII} that $\partial K_0$ is smooth over $V$ and has constant Gaussian curvature equal to $k$ (see also Appendix $A$). In particular, no point of $V$ satisfies the local geodesic condition. This completes the proof.\qed
\medskip
\noindent We thus refine Theorem \procref{TheoremExistence} to obtain:
\proclaim{Lemma \nextprocno}
\noindent Let $\Bbb{H}^{n+1}$ be $(n+1)$-dimensional hyperbolic space, and let $H\subseteq\Bbb{H}^{n+1}$ be a totally geodesic hypersurface. Choose $k>0$, and let $\Omega\subseteq H$ be a bounded open set such that there exists a convex hypersurface $\hat{\Sigma}$ such that:
\medskip
\myitem{(i)} $\hat{\Sigma}$ is a graph below $\Omega$;
\medskip
\myitem{(ii)} the second fundamental form of $\hat{\Sigma}$ is at least $\epsilon$ in the Alexandrov sense, for some $\epsilon>0$; and
\medskip
\myitem{(iii)} the Gaussian curvature of $\Sigma$ is at least $k$ in the Alexandrov sense.
\medskip
\noindent There exists a unique convex, immersed hypersurface $(\Sigma,\partial\Sigma)$ such that:
\medskip
\myitem{(i)} $\Sigma$ is a graph below $\Omega$ and $\partial\Sigma=\partial\Omega$;
\medskip
\myitem{(ii)} $\Sigma$ lies above $\hat{\Sigma}$;
\medskip
\myitem{(iii)} $\Sigma$ has $C^\infty$ interior and is $C^{0,1}$ up to the boundary; and
\medskip
\myitem{(iv)} the Gaussian curvature of $\Sigma$ is equal to $k$.
\medskip
\noindent Moreover if $\partial\Omega$ is smooth, then $\Sigma$ is smooth up to the boundary.
\endproclaim
\proclabel{LemmaDirichletInHyperbolicSpace}
\proof We begin by smoothing the upper barrier. Choose $k'<k$. As in Lemma $2.13$ of \cite{SmiNLD}, there exists a sequence $(\epsilon_n)_\ninn\in ]0,k-k'[$ of positive numbers and a sequence of smooth, convex, immersed hypersurfaces $(\hat{\Sigma}_n)_\ninn$ such that:
\medskip
\myitem{(i)} $(\epsilon_n)_\ninn$ converges to $0$ and $(\hat{\Sigma}_n)_\ninn$ converges to $\hat{\Sigma}$ in the $C^{0,\alpha}$ sense for all $\alpha$;
\medskip
\myitem{(ii)} for all $n$, $\hat{\Sigma}_n$ is a graph over a bounded open subset of $H$; and
\medskip
\myitem{(iii)} for all $n$, the Gaussian curvature of $\hat{\Sigma}_n$ is greater than $k-\epsilon_n$.
\medskip
\noindent Let $(\delta_n)_\ninn>0$ be a sequence of positive numbers converging to $0$. For all $n$, let $H_n$ be the equidistant hypersurface at distance $\delta_n$ from $H$. We may assume that, for all $n$, a portion of $\hat{\Sigma}_n$ is a graph over $H_n$. Let $\Omega_n$ be the subset of $H_n$ over which it as a graph.
\medskip
\noindent For all $n$, since $(\Omega_n,\partial\Omega_n)$ is locally and globally rigid, it follows by Theorem \procref{TheoremExistence} that there exists a smooth, convex hypersurface $\Sigma_n$ which is a graph below $\Omega_n$ such that $\Sigma_n>\hat{\Sigma}_n$ and whose Gaussian curvature is equal to $k'$.
\medskip
\noindent Suppose now that $\partial\Omega$ is smooth. There exists $\epsilon>0$ such that, for all $n$ and for all $P\in\partial\Omega_n$, there exists a geodesic ball $B\subseteq\Omega_n$ such that $P\in\partial B$. For all such $B$, we consider the foliation of constant Gaussian curvature hypersurfaces which are graphs below $B$ and whose boundary is $\partial B$ (in the upper half space model of $\Bbb{H}^{n+1}$, these are merely intersections of spheres in $\Bbb{R}^{n+1}$ with $\Bbb{H}^{n+1}$). Using these foliations and the Geometric Maximum Principle, we find that there exists $\theta>0$ such that, for all $n$, $T\Sigma_n$ makes an angle of at least $\theta$ with $H_n$ along $\partial\Sigma_n$. Bearing in mind the remark following Proposition \procref{PropSecondOrderBoundaryBounds}, this yields uniform lower bounds for the restriction to $\partial\Omega_n$ of the second fundamental form of $\Sigma_n$. Taking limits now yields the desired hypersurface, $\Sigma$.
\medskip
\noindent Consider now the general case. By Lemma \procref{LemmaFirstOrderBounds}, we may nonetheless assume that $(\hat{\Sigma}_n)_\ninn$ converges to a $C^{0,1}$, convex hypersurface $\Sigma$ which is a graph below $\Omega$ such that $\Sigma\geqslant\hat{\Sigma}$. Let $B$ be a geodesic ball such that $\overline{B}\subseteq\Omega$. Using the geometric maximum principle, by considering the foliation of constant Gaussian curvature hypersurfaces which are graphs below $B$ and whose boundary is $B$, we may show that $\Sigma$ lies strictly below $\Omega$ over its interior. We now assert that no point of $\Sigma$ satisfies the local geodesic condition. Indeed, suppose the contrary. By Proposition \procref{PropNodalSetIsClosed} the union of $\partial\Sigma$ with the set of points of $\Sigma\setminus\partial\Sigma$ which satisfy the local geodesic condition is closed. Thus, if $P\in\Sigma$ is such a point, it follows from Lemma \procref{LemmaClosedMeansInConvexHull} that $P$ lies in the convex hull of $\partial\Sigma$. In particular, $P$ lies in $H$ and thus, by convexity, $\Sigma=\Omega$, which is absurd. The assertion follows and it now follows by \cite{CaffII} that $\Sigma$ is smooth over its interior, and this proves existence (see also Appendix $A$).
\medskip
\noindent Let $\Sigma$ be a graph over $\Omega$ of constant Gaussian curvature equal to $k$. Let $f$ be the graph function of $\Sigma$ in conformal coordinates about $H$ (see the example following Proposition \procref{PropFormulaForGaussCurvature}). $f$ satisfies the following equation:
$$
\opDet(f_{;ij} - \opTan(f)(f_{;j}f_{;j} + \delta_{ij}))^{1/n} = k\frac{1}{\opCos(f)^3}(1 + \|\nabla f\|^2)^{(n+2)/2n}.
$$
\noindent Let $\Sigma'$ be another such hypersurface and suppose that $\Sigma'\neq\Sigma$. Let $f'$ be the graph function of $\Sigma'$ in conformal coordinates about $H$. Without loss of generality, there exists $P\in H$ such that $f'(P)>f(P)$ and $f'-f$ is maximised at $P$. Define the field of matrices, $A$, by:
$$
A = (\opHess(f) - \opTan(f)(\nabla f\otimes\nabla f + \opId))^{-1}.
$$
\noindent (This matrix is invertible by convexity of $\Sigma$). $A$ is positive definite. Thus, near $P$, by concavity of $\opDet^{1/n}$, and since $f'>f$, for some $\epsilon,\hat{k}>0$ that we need not calculate:
$$\matrix
\hat{k}\opTr(A^{-1}(f'_{;ij}-f_{;ij})) - \hat{k}\opTan(f)\opTr(A^{-1}(f'_{;i}f'_{;j}-f_{;i}f_{;j}))
\hfill\cr
\qquad\qquad\geqslant \epsilon + \frac{k}{\opCos(f)^3}((1+\|\nabla f'\|^2)^{(n+2)/2} - (1+\|\nabla f\|^2)^{(n+2)/2}).\hfill\cr
\endmatrix$$
\noindent At $P$, since $(f'-f)$ is maximised, $\nabla f'=\nabla f$. Thus, near $P$;
$$
\opTr(A^{-1}(f'_{;ij}-f_{;ij})) > 0
$$
\noindent This yields a contradiction by the Maximum Principle. Uniqueness follows and this completes the proof.\qed
\goodbreak
\newhead{Relations to Existing Results}
\noindent With small modifications, these techniques may be adapted to yield existing results. First, considering $\Bbb{R}^n$ as a subspace of $\Bbb{R}^{n+1}$ in the natural manner, we recover the following theorem of Guan (see \cite{Guan}), which is the analogue in Euclidean space of Lemma \procref{LemmaDirichletInHyperbolicSpace}:
\headlabel{HeadRelationsToExistingResults}
\proclaim{Theorem \nextprocno, {\bf [Guan, 1998]}}
\noindent Choose $k>0$, and let $\Omega\subseteq\Bbb{R}^n$ be a bounded open set. Suppose that there exists a convex hypersurface, $\hat{\Sigma}$ such that:
\medskip
\myitem{(i)} $\hat{\Sigma}$ is a graph below $\Omega$;
\medskip
\myitem{(ii)} the second fundamental form of $\hat{\Sigma}$ is at least $\epsilon$ in the Alexandrov sense, for some $\epsilon>0$; and
\medskip
\myitem{(iii)} the Gaussian curvature of $\Sigma$ is at least $k$ in the Alexandrov sense.
\medskip
\noindent There exists a unique convex, immersed hypersurface $(\Sigma,\partial\Sigma)$ such that:
\medskip
\myitem{(i)} $\Sigma$ is a graph below $\Omega$ and $\partial\Sigma=\partial\Omega$;
\medskip
\myitem{(ii)} $\Sigma$ lies above $\hat{\Sigma}$;
\medskip
\myitem{(iii)} $\Sigma$ has $C^\infty$ interior and is $C^{0,1}$ up to the boundary; and
\medskip
\myitem{(iv)} the Gaussian curvature of $\Sigma$ is equal to $k$.
\medskip
\noindent Moreover, if $\partial\Omega$ is smooth, then $\Sigma$ is smooth up to the boundary.

\endproclaim
\proclabel{TheoremGuan}
\remark Although, as in Lemma \procref{LemmaDirichletInHyperbolicSpace}, if we identify $(\Sigma_0,\partial\Sigma_0)=(\Omega,\partial\Omega)$, then the Gauss Curvature Equation is not elliptic at $f_0=0$, this, in itself, does not present a serious difficulty since there exist functions arbitrarily close to $f_0$ where the Gauss Curvature Equation is elliptic. The particular difficulty in Euclidean space lies in obtaining functions near $f_0$ for which the Gauss Curvature Equation is also stable. We circumvent this by approximating $\Bbb{R}^n$ by spaces of constant negative sectional curvature.
\medskip
\proof Using polar coordinates for $\Bbb{R}^n$, we identify $\Bbb{R}^{n+1}$ with $\Sigma^{n-1}\times]0,\infty[\times\Bbb{R}$, where $\Sigma^{n+1}$ is the unit sphere. We thus denote a point in $\Bbb{R}^{n+1}$ by the coordinates $(\theta,r,t)\in\Sigma^{n-1}\times]0,\infty[\times\Bbb{R}$. Let $g^\Sigma$ denote the standard metric over $\Sigma^{n-1}$. For $\epsilon>0$, we define the metric $g_\epsilon$ over $\Bbb{R}^{n+1}$ such that, at $(\theta,r,t)$:
$$
g = \opCosh^2(\epsilon t)( \opSinh^2(\epsilon r)g^{\Sigma} \oplus dr^2) \oplus dt^2.
$$
\noindent This metric is smooth and has constant curvature equal to $-\epsilon$. Indeed, this formula is obtained by using polar coordinates for $\Bbb{H}^n$ about a point and subsequently by identifying $\Bbb{H}^{n+1}$ with $\Bbb{H}^n\times\Bbb{R}$ using the foliation by geodesics normal to a totally geodesic hypersurface.
\medskip
\noindent With respect to this metric, $\Bbb{R}^n$ is identified with a totally geodesic hypersurface, and, for all $k'<k$, there exists $\epsilon>0$ such that $\hat{\Sigma}$ satisfies the hypotheses of Lemma \procref{LemmaDirichletInHyperbolicSpace}, with $k'$ instead of $k$. There therefore exists $\Sigma_\epsilon\subseteq\Bbb{R}^{n+1}$ possessing the desired properties and of constant Gaussian curvature equal to $k'$ with respect to the metric $g_\epsilon$. Existence follows by taking limits as in Lemma \procref{LemmaDirichletInHyperbolicSpace}.
\medskip
\noindent To prove uniqueness, let $\Sigma_1$ and $\Sigma_2$ be two solutions. Suppose that $\Sigma_1\neq\Sigma_2$. Without loss of generality, there is a point of $\Sigma_1$ lying below $\Sigma_2$. There therefore exists a translate $\Sigma_1'$ of $\Sigma_1$ in the vertical direction which lies strictly above $\Sigma_1$ and which is an exterior tangent to $\Sigma_2$ at some point $P'$. Since $\partial\Sigma_1'$ lies strictly above $\partial\Sigma_2$, $P'$ is an interior point of $\Sigma_1'$ and $\Sigma_2$. It follows by the strong Geometric Maximum Principle that $\Sigma_1'$ and $\Sigma_2$ coincide, which is absurd. Uniqueness follows and this completes the proof.\qed
\medskip
\noindent If $M$ is a Riemannian manifold, we say that a bounded open subset $\Omega\subseteq M$ satisfies a uniform exterior ball condition if and only if there exists $\epsilon>0$ such that for every $P\in\partial\Omega$, there exists an open geodesic ball $B\subseteq\Omega^c$ of radius $\epsilon$ such that:
$$
P\in \partial B\minter\partial\Omega.
$$
\noindent By compactness, $\Omega$ satisfies a uniform exterior ball condition for a given metric over $M$ if and only if it satisfies this condition for any metric over $M$, and we thus extend this condition to subsets of arbitrary $C^\infty$ manifolds.
\medskip
{\sl\noindent Example:\ }Any compact, open subset with smooth boundary satisfies a uniform exterior ball condition.\qed
\medskip
{\sl\noindent Example:\ }Any convex, open subset satisfies a uniform exterior ball condition.\qed
\medskip
\noindent We now recover the following theorem of Rosenberg and Spruck (see \cite{RosSpruck}), which has also recently been proven in a more general setting by Guan, Spruck and Szapiel (see \cite{GuanSpruckII}):
\proclaim{Theorem \nextprocno, {\bf [Rosenberg, Spruck, (1994)]}}
\noindent Let $\Omega\subseteq\partial_\infty\Bbb{H}^{n+1}$ be a non-trivial open subset whose boundary satisfies the uniform exterior ball condition. Then, for all $k\in]0,1[$, there exists a convex, immersed hypersurface $\Sigma_k\subseteq\Bbb{H}^{n+1}$ such that:
\medskip
\myitem{(i)} identifying $\partial_\infty\Bbb{H}^{n+1}$ with $\Bbb{R}^n\munion\left\{\infty\right\}$ and viewing $\Omega$ as a subset of $\Bbb{R}^n$, $\Sigma_k$ is a graph over $\Omega$;
\medskip
\myitem{(ii)} $\Sigma_k$ is smooth and $C^{0,1}$ up to the boundary;
\medskip
\myitem{(iii)} $\partial\Sigma_k=\partial\Omega$; and
\medskip
\myitem{(iv)} $\Sigma_k$ has constant Gaussian curvature equal to $k$.
\medskip
\noindent Moreover, if $\Omega$ is star-shaped, then $\Sigma_k$ is unique.
\endproclaim
\proclabel{TheoremRosSpruck}
\remark In this case, we use horospheres as upper barriers. Since these have curvature equal to $1$, we can only prove existence for hypersurfaces of curvature less than $1$, hence the hypothesis on $k$.
\medskip
\proof We identify $\Bbb{H}^{n+1}$ with the upper half space $\Bbb{R}^n\times]0,\infty[$ in the standard manner. We thus identify $\partial_\infty\Bbb{H}^{n+1}$ with $\Bbb{R}^n\munion\left\{\infty\right\}$ and view $\Omega$ as a subset of $\Bbb{R}^n$. For $\epsilon>0$, let $H_\epsilon=\Bbb{R}^n\times\left\{\epsilon\right\}$ be the horosphere at height $\epsilon$ above $\Bbb{R}^n$. We define $\Omega_\epsilon\subseteq H_\epsilon$ by
$$
\Omega_\epsilon = \left\{(x,\epsilon)\text{ s.t. }x\in\Omega\right\}.
$$
\noindent By the uniform exterior ball condition, for $\epsilon$ sufficiently small, $\partial\Omega_\epsilon$ is uniformly strictly convex as a subset of $\Bbb{H}^{n+1}$ with respect to the outward pointing unit normal in $H_\epsilon$.
\medskip
\noindent Let $K_\epsilon$ be the complement of $\Omega_\epsilon$ in $H_\epsilon$. Let $\hat{K}_\epsilon$ be the convex hull of $K_\epsilon$ in $\Bbb{H}^{n+1}$. We denote by $\Sigma_{0,\epsilon}$ the portion of $\partial\hat{K}_\epsilon$ lying above $H_\epsilon$. In other words:
$$
\partial\hat{K}_\epsilon = (\partial\hat{K}_\epsilon\minter H_\epsilon)\munion\Sigma_{0,\epsilon}.
$$
\noindent Since it is locally ruled, $\Sigma_{0,\epsilon}$ serves as a lower barrier for the problem (see \cite{SmiFCS}). We define $(\hat{\Sigma}_\epsilon,\partial\hat{\Sigma}_\epsilon)=(\Omega,\partial\Omega)$. The only difference between our current framework and that of Theorem \procref{TheoremExistence} is that it is the upper barrier, $\hat{\Sigma}_\epsilon$ that is smooth and the lower barrier, $\Sigma_{0,\epsilon}$ that is not. The only change required to adapt the proof to our framework is therefore to replace $(f-f_0)$ in Corollary \procref{CorSecondBarrierEstimate} with $(f-\hat{f})$. The uniform strict convexity of $\Omega_\epsilon$ as a subset of $\Bbb{H}^{n+1}$ with respect to the normal in $H_\epsilon$ ensures uniform lower bounds of the restriction to $\partial\Omega$ of the second fundamental form of any surface $\Sigma$ which is a graph above $\Omega$ such that $\partial\Sigma=\partial\Omega$. Thus proceeding as in Theorem \procref{TheoremExistence}, we show that there exists a graph $\Sigma_\epsilon$ over $\Omega_\epsilon$ which is smooth up to the boundary and has constant Gaussian curvature equal to $k$.
\medskip
\noindent Taking limits yields a $C^{0,1}$ graph $\Sigma$ over $\Omega$ such that $\partial\Sigma=\partial\Omega$. Let $Y\subseteq\Sigma$ be the set of all points satisfying the local geodesic condition. By Proposition \procref{PropNodalSetIsClosed}, $\partial\Sigma\munion Y$ is closed. It follows as in Lemma \procref{LemmaClosedMeansInConvexHull} that $Y$ is contained in the convex hull of the intersection of some totally geodesic hyperplane $H$ with $\partial\Omega$ (see \cite{SmiNLD} for details). In particular, if $Y$ is non-empty, then, viewed as a graph, $\Sigma$ is vertical at some point on the boundary. However, consider a point $P\in\partial\Omega_\epsilon$ and a geodesic ball $B\subseteq H_\epsilon$ such that $B\subseteq\Omega^c$ and $P\in\partial B$. Using the foliation of constant Gaussian curvature hypersurfaces in $\Bbb{H}^n$ whose boundary coincides with $\partial B$, we deduce by the Geometric Maximum Principle that there exists $\theta>0$ such that, for $\epsilon$ sufficiently small, $\Sigma_\epsilon$ makes an angle at $P$ of at least $\theta$ with the foliation of vertical geodesics along $\partial\Omega$. Moreover, $\theta$ may be chosen independant of $P$. Taking limits, it follows that $\Sigma$ is everywhere strictly convex and is therefore smooth over the interior by \cite{CaffII}. This proves existence.
\medskip
\noindent Suppose now that $\Omega$ is star-shaped, and let $\Sigma_1$ and $\Sigma_2$ be two solutions. Suppose that $\Sigma_1\neq\Sigma_2$. Without loss of generality, there exists a point $P\in\Sigma_1$ lying below $\Sigma_2$. As before, we identify $\Bbb{H}^{n+1}$ with $\Bbb{R}^n\times]0,\infty[$. Without loss of generality, we may suppose that $\Omega$ is star-shaped about $(0,0)$. Consider the family $(M_\lambda)_{\lambda>1}$ of isometries of $\Bbb{H}^{n+1}$ given by:
$$
M_\lambda(x,t) = (\lambda x,\lambda t).
$$
\noindent There exists $\lambda>1$ such that $M_\lambda \Sigma_1$ is an exterior tangent to $\Sigma_2$ at some point $P'$. Since $M_\lambda\partial\Sigma_1\minter\partial\Sigma_2=\emptyset$, $P'$ is an interior point of $\Sigma_1$ and $\Sigma_2$. It follows by the strong Geometric Maximum Principle that $M_\lambda\Sigma_1=\Sigma_2$, which is absurd. Uniqueness follows and this completes the proof.\qed
\inappendicestrue
\global\headno=0
\goodbreak
\newhead{Regularity of Limit Hypersurfaces}
\noindent Let $M^{n+1}$ be an $(n+1)$-dimensional Riemannian manifold. Choose $k>0$ let $(\Sigma_m)_\minn$ be a sequence of smooth, convex, immersed hypersurfaces in $M$ of constant Gaussian curvature equal to $k$. Suppose that there exists a $C^{0,1}$ locally convex, immersed hypersurface, $\Sigma_0$ to which $(\Sigma_m)_\minn$ converges in the $C^{0,\alpha}$ sense for all $\alpha$. For all $m\in\Bbb{N}$, let $\msf{N}_m$ and $A_m$ be the unit normal vector field and the second fundamental form respectively of $\Sigma_m$. Choose $p_0\in\Sigma_0$ and let $(p_m)_\minn\in(\Sigma_m)_\minn$ be a sequence converging to $p_0$. For all $r>0$, and for all $m\in\Bbb{N}\munion\left\{0\right\}$, let $B_{m,r}$ be the ball of radius $r$ (with respect to the intrinsic metric) about $p_m$ in $\Sigma_m$.
\medskip
\noindent We will say that $\Sigma_0$ is functionally strictly convex at $p_0$ if there exists a smooth function, $\varphi$, defined on $M$ near $p_0$ such that:
\medskip
\myitem{(i)} $\varphi$ is strictly convex;
\medskip
\myitem{(ii)} $\varphi(p_0)>0$; and
\medskip
\myitem{(iii)} the connected component of $\varphi^{-1}([0,\infty[)\minter\Sigma_0$ containing $p_0$ is compact.
\medskip
\noindent Observe that if $M$ is a space form, then $\Sigma_0$ is functionally strictly convex whenever it is strictly convex. We will prove:
\proclaim{Lemma \nextprocno}
\noindent If $\Sigma_0$ is funtionally strictly convex at $p_0$, then there exists $r>0$ such that $(B_{m,r},p_m)_\minn$ converges to $(B_{0,r},p_0)$ in the pointed $C^\infty$-Cheeger Gromov sense. In particular, $B_{0,r}$ is a smooth, convex immersion of constant Gaussian curvature equal to $k$.
\endproclaim
\proclabel{LemmaSmoothnessOfLimit}
\noindent As in section \headref{HeadSecondOrderBoundaryBounds}, we denote by $\Cal{B}$ the family of constants which depend continuously on the data: $M$, $k$, $(\Sigma_0,p_0)$ and the $C^1$ jets of $(\Sigma_m,p_m)_\minn$. In this section, for any positive quantity, $X$, we denote by $O(X)$ any term which is bounded in magnitude by $K\left|X\right|$ for some $K$ in $\Cal{B}$.
\medskip
\noindent The following elementary lemma will be of use in the proof:
\proclaim{Lemma \nextprocno}
\noindent For $\lambda>0$ and for all $a,b\in\Bbb{R}$:
$$
(a+b)^2 \leqslant (1+\lambda)a^2 + (1+\lambda^{-1})b^2.
$$
\endproclaim
\proclabel{LemmaBasicSquareRelation}
{\bf\noindent Proof of Lemma \procref{LemmaSmoothnessOfLimit}:\ }Since $(\Sigma_m)_\ninn$ converges to $\Sigma_0$ and since $\Sigma_0$ is functionally strictly convex at $p_0$, there exists $\epsilon,h>0$, open sets $\Omega_0,(\Omega_m)_\minn\subseteq M$ and, for every $m$, a smooth function $\varphi_m:\Omega_m\rightarrow[0,\infty[$ such that:
\medskip
\myitem{(i)} for all $m$, $\Omega_m$ is a neighbourhood of $p_m$ and $(\Omega_m)_\minn$ converges to $\Omega_0$ in the Hausdorff sense;
\medskip
\myitem{(ii)} $(\varphi_m)_\minn$ converges to $\varphi_0$ in the $C^\infty$ sense;
\medskip
\myitem{(iii)} for all $m$, $\opHess(\varphi_m)\geqslant\epsilon\opId$;
\medskip
\myitem{(iii)} for all $m$, $\varphi_m(p_0)=2h$; and
\medskip
\myitem{(iv)} for all $m$, the connected component of $p_m$ in $\Sigma_m\minter\Omega_m$ is compact with smooth boundary and $\varphi_m$ equals zero along the boundary: we denote this connected component by $\Sigma_{m,0}$.
\medskip
\noindent We may assume that, for all $m$, $\varphi_m\leqslant 1$ over $\Sigma_{m,0}$. Finally, after reducing $\epsilon$ if necessary, there exists a smooth, unit length vector field $X$ defined over a neighbourhood of $p_0$ such that, for all $m$, throughout $\Sigma_{m,0}$, $\langle X,\msf{N}_m\rangle\geqslant\epsilon$. We now follow an adaptation of reasoning presented by Pogorelov in \cite{Pog}.
\medskip
\noindent Choose $\alpha\geqslant 1$. For all $m$, we define the function $\Phi_m$ by:
$$
\Phi_m = \alpha\opLog(\varphi_m) - \langle X,\msf{N}_m\rangle + \opLog(\|A_m\|),
$$
\noindent where $\|A_m\|$ is the operator norm of $A_m$, which is equal to the highest eigenvalue of $A_m$. We aim to obtain a priori upper bounds for $\Phi_m$ for some $\alpha$. We trivially obtain a-priori bounds whenever $\|A_m\|\leqslant 1$. We thus consider the region where $\|A_m\|\geqslant 1$. Choose $m\in\Bbb{N}$ and $P\in\Sigma_{m,0}$. Let $\lambda_1\geqslant...\geqslant\lambda_n$ be the eigenvalues of $A_m$ at $P$. In particular, $\lambda_1=\|A_m\|$. Let $e_1,...,e_n$ be the corresponding orthonormal basis of eigenvectors. In the sequel, we will suppress $m$.
\medskip
\noindent Let the subscript ``$;$'' denote covariant differentiation with respect to the Levi-Civita covariant derivative of $\Sigma$. Thus, for example:
$$
A_{ij;k} = (\nabla^\Sigma_{e_k}A)(e_i,e_j).
$$
\noindent We consider the Laplacian, $\Delta$, defined on functions by:
$$
\Delta f = \sum_{i=1}^n\frac{1}{\lambda_i}f_{;ii}.
$$
\noindent We aim to use the Maximum Principle in conjunction with $\Delta$. Thus, in the sequel, we will only be interested in the orders of magnitude of potentially negative terms.
\medskip
\noindent In analogy to Corollary \procref{CorSuperharmonic}, at $P$:
$$
\Delta\opLog(\lambda_1) \geqslant \sum_{i,j=1}^n\frac{1}{\lambda_1\lambda_i\lambda_j}A_{ij;1}A_{ij;1} - \sum_{i=1}^n\frac{1}{\lambda_1\lambda_1\lambda_i}A_{11;i}A_{11;i} - O(1) - O(\|A^{-1}\|),
$$
\noindent in the weak sense. However, by Lemma \procref{LemmaCommutationRelations}, for all $i$:
$$
A_{11;i} = A_{i1;1} + R^M_{i1\nu1},
$$
\noindent where $\nu$ represents the exterior normal direction to $\Sigma$. Thus, bearing in mind Lemma \procref{LemmaBasicSquareRelation}, and that $\lambda_1\geqslant 1$, we obtain:
$$
\Delta\opLog(\lambda_1) \geqslant \sum_{i=2}^n\frac{1}{2\lambda_1\lambda_1\lambda_i}A_{i1;1}A_{i1;1} - O(1) - O(\|A^{-1}\|),
$$
\noindent in the weak sense. Differentiating the Gauss Curvature Equation yields, for all $j$:
$$
\sum_{i=1}^n\frac{1}{\lambda_i}A_{ii;j} = 0,
$$
\noindent Thus:
$$\matrix
-\Delta\langle X,\msf{N}\rangle \hfill&\geqslant \langle X,\msf{N}\rangle\opTr(A) - O(1) - O(\|A^{-1}\|)\hfill\cr
&\geqslant \epsilon\lambda_1 - O(1) - O(\|A^{-1}\|).\hfill\cr
\endmatrix$$
\noindent Finally:
$$
\Delta(\alpha\opLog(\varphi)) \geqslant \frac{\alpha}{\varphi}\epsilon\opTr(A^{-1}) - \frac{\alpha}{\varphi^2}\sum_{i=1}^n\frac{1}{\lambda_i}\varphi_{;i}\varphi_{;i} - O(\alpha).
$$
\noindent However, bearing in mind Lemma \procref{LemmaCommutationRelations}:
$$
\Phi_{;i} = \frac{\alpha}{\varphi}\varphi_{;i} - {X^\nu}_{;i} - X^i\lambda_i + \frac{1}{\lambda_1}A_{i1;1} + \frac{1}{\lambda_1}R^M_{i1\nu1},
$$
\noindent where $\nu$ is the exterior normal direction over $\Sigma$. Thus, by induction on Lemma \procref{LemmaBasicSquareRelation}, modulo $\nabla\Phi$:
$$
\left|\frac{\alpha}{\varphi}\varphi_{;i}\right|^2 \leqslant \frac{4}{\lambda_1^2}A_{i1;1}A_{i1;1} + \frac{4}{\lambda_1^2}(R^M_{i1\nu 1})^2 + 4(X^i\lambda_i)^2 + 4(X^\nu_{;i})^2.
$$
\noindent Thus, bearing in mind that $\lambda_1\geqslant\lambda_i$ for all $i$ and that $\lambda_1\geqslant 1$, we obtain, modulo $\nabla\Phi$:
$$
\frac{\alpha}{\varphi^2}\sum_{i=2}^n\frac{1}{\lambda_i}\varphi_{;i}\varphi_{;i}
=O(\alpha^{-1}\|A^{-1}\|) + O(\alpha^{-1}\lambda_1) + \sum_{i=2}^n\frac{4}{\alpha\lambda_1^2\lambda_i}A_{i1;1}A_{i1;1}.
$$
\noindent Since $\varphi$ is bounded above (and thus $\varphi^{-1}$ is bounded below), for sufficently large $\alpha$ we obtain, modulo $\nabla\Phi$:
$$\matrix
\Delta\Phi \hfill&\geqslant \frac{\epsilon}{2}\lambda_1 - O(\lambda_1^{-1}\varphi^{-2}) - O(1)\hfill\cr
&= (\varphi^{2\alpha}\|A\|)^{-1}(\frac{\epsilon}{2}(\varphi^\alpha\|A\|)^2 - O(\varphi^\alpha\|A\|) - O(1)).\hfill\cr
\endmatrix$$
\noindent There therefore exists $K_1>0$ in $\Cal{B}$ such that if $(\varphi^\alpha\|A\|)\geqslant K$, then the right hand side is positive. However, for all $m\in\Bbb{N}$, $\Phi_m=-\infty$ along $\partial\Sigma_{m,0}$. There thus exists a point $P\in\Sigma_{m,0}$ where $\Phi_m$ is maximised. By the Maximum Principle, at this point, either $\|A\|\leqslant 1$ or $\varphi^\alpha\|A\|\leqslant K_1$. Taking exponentials, there therefore exists $K_2>0$ in $\Cal{B}$ such that, for all $m\in\Bbb{N}$, throughout $\Sigma_{m,0}$:
$$
\varphi^\alpha\langle X,\msf{N}_m\rangle^{-1}\|A_m\|\leqslant K_2.
$$
\noindent Since $\langle X,\msf{N}_m\rangle\leqslant 1$, this yields a-priori bounds for $\|A_m\|$ over the intersection of $\Sigma_{m,0}$ with $\varphi_m\geqslant h$. Using, for example, an adaptation of the proof of Theorem $1.2$ of \cite{SmiSLC} in conjunction with the Bernstein Theorem \cite{CalI}, \cite{Jorgens} \& \cite{Pog} of Calabi, J\"orgens, Pogorelov, we obtain a-priori $C^k$ bounds for $\Sigma_{m,0}$ over the region $\varphi_m\geqslant 3h$ for all $k$. The result now follows by the Arzela-Ascoli Theorem.\qed
\goodbreak
\newhead{Bibliography}
{\leftskip = 5ex \parindent = -5ex
\leavevmode\hbox to 4ex{\hfil \cite{CaffI}}\hskip 1ex{Caffarelli L., A localization property of viscosity solutions to the Monge-Amp\`ere equation and their strict convexity, {\sl Annals of Math.} {\bf 131} (1990), 129---134}
\medskip
\leavevmode\hbox to 4ex{\hfil \cite{CaffII}}\hskip 1ex{Caffarelli L., {\sl Monge-Amp\'ere equation, div-curl theorems in Lagrangian coordinates, compression and rotation}, Lecture Notes, 1997}
\medskip
\leavevmode\hbox to 4ex{\hfil \cite{CaffNirSprI}}\hskip 1ex{Caffarelli L., Nirenberg L., Spruck J., The Dirichlet problem for nonlinear second-Order elliptic equations. I. Monge Amp\`ere equation, {\sl Comm. Pure Appl. Math.} {\bf 37} (1984), no. 3, 369--402}
\medskip
\leavevmode\hbox to 4ex{\hfil \cite{CaffNirSprII}}\hskip 1ex{Caffarelli L., Kohn J. J., Nirenberg L., Spruck J., The Dirichlet problem for nonlinear second-order elliptic equations. II. Complex Monge Amp\`ere, and uniformly elliptic, equations, {\sl Comm. Pure Appl. Math.} {\bf 38} (1985), no. 2, 209--252}
\medskip
\leavevmode\hbox to 4ex{\hfil \cite{CaffNirSprV}}\hskip 1ex{Caffarelli L., Nirenberg L., Spruck J., Nonlinear second-order elliptic equations. V. The Dirichlet problem for Weingarten hypersurfaces, {\sl Comm. Pure Appl. Math.} {\bf 41} (1988), no. 1, 47--70}
\medskip
\leavevmode\hbox to 4ex{\hfil \cite{CalI}}\hskip 1ex{Calabi E., Improper affine hypersurfaces of convex type and a generalization of a theorem by K. J\"orgens, {\sl Mich. Math. J.} {\bf 5} (1958), 105--126}
\medskip
\leavevmode\hbox to 4ex{\hfil \cite{Gauss}}\hskip 1ex{Gauss C. F., General investigations of curved surfaces, Translated from the Latin and German by Adam Hiltebeitel and James Morehead, Raven Press, Hewlett, N.Y.}
\medskip
\leavevmode\hbox to 4ex{\hfil \cite{GilbTrud}}\hskip 1ex{Gilbarg D., Trudinger N. S., {\sl Elliptic Partial Differential Equations of Second Order}, Classics in Mathematics, Springer, Berlin, (2001)}
\medskip
\leavevmode\hbox to 4ex{\hfil \cite{Guan}}\hskip 1ex{Guan B., The Dirichlet problem for Monge-Amp\`ere equations in non-convex domains and spacelike hypersurfaces of constant Gauss curvature, {\sl Trans. Amer. Math. Soc.} {\bf 350} (1998), 4955--4971}
\medskip
\leavevmode\hbox to 4ex{\hfil \cite{GuanSpruckO}}\hskip 1ex{Guan B., Spruck J., Boundary value problems on $\Bbb{S}^n$ for surfaces of constant Gauss curvature, {\sl Ann. of Math.} {\bf 138} (1993), 601--624}
\medskip
\leavevmode\hbox to 4ex{\hfil \cite{GuanSpruckI}}\hskip 1ex{Guan B., Spruck J., The existence of hypersurfaces of constant Gauss curvature with prescribed boundary, {\sl J. Differential Geom.} {\bf 62} (2002), no. 2, 259--287}
\medskip
\leavevmode\hbox to 4ex{\hfil \cite{GuanSpruckII}}\hskip 1ex{Guan B., Spruck J., Szapiel M., Hypersurfaces of constant curvature in Hyperbolic space I, {\sl J. Geom. Anal}}
\medskip
\leavevmode\hbox to 4ex{\hfil \cite{GuanSpruckIII}}\hskip 1ex{Guan B., Spruck J., Hypersurfaces of constant curvature in Hyperbolic space II, arXiv:0810.1781}
\medskip
\leavevmode\hbox to 4ex{\hfil \cite{Gut}}\hskip 1ex{Guti\'errez C., {\sl The Monge-Amp\`ere equation}, Progress in Nonlinear Differential Equations and Their Applications, {\bf 44}, Birkh\"a user, Boston, (2001)}
\medskip
\leavevmode\hbox to 4ex{\hfil \cite{Jorgens}}\hskip 1ex{J\"orgens K., \"Uber die L\"osungen der Differentialgleichung $rt-s^2=1$, {\sl Math. Ann.} {\bf 127} (1954), 130--134}
\medskip
\leavevmode\hbox to 4ex{\hfil \cite{LabI}}\hskip 1ex{Labourie F., Un lemme de Morse pour les surfaces convexes (French), {\sl Invent. Math.} {\bf 141} (2000), no. 2, 239--297}
\medskip
\leavevmode\hbox to 4ex{\hfil \cite{LoftI}}\hskip 1ex{Loftin J. C., Affine spheres and convex $\Bbb{RP}\sp n$-manifolds, {\sl Amer. J. Math.} {\bf 123} (2001), no. 2, 255--274}
\medskip
\leavevmode\hbox to 4ex{\hfil \cite{LoftII}}\hskip 1ex{Loftin J. C., Riemannian metrics on locally projectively flat manifolds, {\sl Amer. J. Math.} {\bf 124} (2002), no. 3, 595--609}
\medskip
\leavevmode\hbox to 4ex{\hfil \cite{Pog}}\hskip 1ex{Pogorelov A. V., On the improper affine hyperspheres, {\sl Geom. Dedicata} {\bf 1} (1972), 33--46}
\medskip
\leavevmode\hbox to 4ex{\hfil \cite{PogB}}\hskip 1ex{Pogorelov A. V., Extrinsic geometry of convex surfaces, Translated from the Russian by Israel Program for Scientific}
\medskip
\leavevmode\hbox to 4ex{\hfil \cite{RosSpruck}}\hskip 1ex{Rosenberg H., Spruck J., On the existence of convex hypersurfaces of constant Gauss curvature in hyperbolic space, {\sl J. Differential Geom.} {\bf 40} (1994), no. 2, 379--409}
\medskip
\leavevmode\hbox to 4ex{\hfil \cite{SmiPPH}}\hskip 1ex{Smith G., Compactness results for immersions of prescribed Gaussian curvature II - geometric aspects, arXiv:1002.2982}
\medskip
\leavevmode\hbox to 4ex{\hfil \cite{SmiPPG}}\hskip 1ex{Smith G., The Plateau problem for general curvature functions, arXiv:1008.3545}
\medskip
\leavevmode\hbox to 4ex{\hfil \cite{SmiNLD}}\hskip 1ex{Smith G., The Non-Linear Dirichlet Problem in Hadamard Manifolds,\break arXiv:0908.3590}
\medskip
\leavevmode\hbox to 4ex{\hfil \cite{SmiFCS}}\hskip 1ex{Smith G., Moduli of Flat Conformal Structures of Hyperbolic Type, arXiv:0804.0744}
\medskip
\leavevmode\hbox to 4ex{\hfil \cite{SmiSLC}}\hskip 1ex{Smith G., Special Lagrangian curvature, arXiv:math/0506230}
\medskip
\leavevmode\hbox to 4ex{\hfil \cite{Spruck}}\hskip 1ex{Spruck J., Fully nonlinear elliptic equations and applications to geometry, {\sl Proceedings of the International Congress of Mathematicians}, Vol. 1, 2 (Z\"urich, 1994), 1145--1152, Birkh\"auser, Basel, (1995)}
\medskip
\leavevmode\hbox to 4ex{\hfil \cite{TrudWang}}\hskip 1ex{Trudinger N. S., Wang X., On locally convex hypersurfaces with boundary, {\sl J. Reine Angew. Math.} {\bf 551} (2002), 11--32}
\par}

\enddocument